\documentclass{amsart}
\usepackage{lmodern}
\usepackage[utf8]{inputenc}
\usepackage[T1]{fontenc}
\usepackage[english]{babel}
\usepackage{amsmath,amssymb,amsthm}
\usepackage{mathtools} 
\usepackage{bm} 
\usepackage{esint} 
\usepackage{dsfont} 
\usepackage{color}
\usepackage{csquotes}
\usepackage{enumitem}
\usepackage[colorlinks=true]{hyperref}
\usepackage[letterpaper,nomarginpar]{geometry}
\usepackage{pgf,tikz}
\usetikzlibrary{arrows}
\usepackage[foot]{amsaddr}

\DeclareMathOperator{\Dom}{Dom}
\DeclareMathOperator{\Leb}{Leb}
\DeclarePairedDelimiter{\abs}{\lvert}{\rvert}
\DeclarePairedDelimiter{\brc}{(}{)}
\DeclarePairedDelimiter{\sqbrc}{[}{]}
\DeclarePairedDelimiter{\norm}{\lVert}{\rVert}

\DeclareMathOperator{\vect}{Vect}

\DeclareMathOperator{\Cu}{curl}

\newcommand{\interval}[4]{\mathopen{#1}#2\mathclose{},#3\mathclose{#4}}
\newcommand{\intcc}[2]{\interval{[}{#1}{#2}{]}}
\newcommand{\intoc}[2]{\interval{]}{#1}{#2}{]}}

\newcommand{\numberset}[1]{\mathbb{#1}}
\newcommand{\N}{\numberset{N}}
\newcommand{\R}{\numberset{R}}
\newcommand{\Sp}{\numberset{S}}

\newcommand{\cC}{\mathcal{C}}
\newcommand{\cM}{\mathcal{M}}

\newcommand{\diff}{\mathop{}\mathopen{}\mathrm{d}} 
\newcommand{\indic}{\mathds{1}}
\newcommand{\ie}{i.e.\ }

\newcommand{\hs}{\mathcal{H}}
\newcommand{\mres}{\mathop{\hbox{\vrule height 6pt width .5pt depth 0pt \vrule height .5pt width 4pt depth 0pt}}\nolimits}
\newcommand{\ve}{\varepsilon}

\numberwithin{equation}{section}

\theoremstyle{plain}
\newtheorem{thm}{Theorem}[section]
\newtheorem{lem}[thm]{Lemma}
\newtheorem{prop}[thm]{Proposition}
\newtheorem{cor}[thm]{Corollary}

\theoremstyle{definition}
\newtheorem{df}{Definition}[section]

\theoremstyle{remark}
\newtheorem{rmk}{Remark}[section]

\begin{document}
\author{Pierre Bochard$^\dagger$}
\address{$^\dagger$Université Lyon 1, CNRS UMR 5208, Institut Camille Jordan, 69622 Villeurbanne, France}
\email{bochard@math.univ-lyon1.fr}
\author{Paul Pegon$^\ddagger$}
\address{$^\ddagger$Laboratoire de Mathématiques d'Orsay, Univ. Paris-Sud, CNRS, Université Paris-Saclay, 91405 Orsay Cedex, France}
\email{paul.pegon@math.u-psud.fr}

\title[Regularity of curl-free vector fields]{A kinetic selection principle for curl-free vector
fields of unit norm}

\begin{abstract}
This article is devoted to the generalization of results obtained in 2002 in 
\cite{Jabin_Otto_Perthame_Line_energy_2002} by Jabin, Otto and Perthame. In their article they proved that planar vector fields taking value into the unit sphere of the euclidean norm and satisfying a given kinetic equation are locally Lipschitz. Here, we study the same question replacing the unit sphere of the euclidean norm by the unit sphere of \emph{any} norm. Under natural asumptions on the norm, namely smoothness and a qualitative convexity property, that is to be of power type $p$, we prove that planar vector fields taking value into the unit sphere of such a norm and satisfying a certain kinetic equation are locally $\frac{1}{p-1}$-H\"older continuous. Furthermore we completely describe the behaviour of such a vector field around singular points as a \emph{vortex} associated to the norm. As our kinetic equation implies for the vector field to be curl-free, this can be seen as a selection principle for curl-free vector fields valued in spheres of general norms which rules out line-like singularities.
\end{abstract}

\maketitle
\tableofcontents
\section*{Preliminaries and notations}
Our notations in this article will be standard. Beware that we will consider in this article maps instead of equivalence classes for the almost everywhere relation; all maps will be taken Borel 
measurable. More specifically, let $\Omega, E \subset \R^2$. We will denote by $\mathcal{L}^1_{\mathrm{loc}}(\Omega,E)$ the set of Borel maps $m$ which are locally integrable and such that $m(x) \in B$ for all $x$ in $\Omega$. Of course we do not lose generality with this choice. For $x\in \R^2$, $x^\perp$ will denote its rotation by an angle $+\pi/2$. The euclidean norm will be denoted by $\abs{\cdot}$ and a general norm by $\norm{\cdot}$.

\section{Introduction}
The goal of this paper is the generalization of regularity results on solutions of 
a kinetic equation
obtained by Jabin, Otto and Perthame in their article 
\cite{Jabin_Otto_Perthame_Line_energy_2002}. In order for the reader to understand where we are
starting from, we will begin by briefly recalling their results.
Let $\Omega \subset \R^2$ be an open set. Jabin, Otto and Perthame were interested in the study of 
the following micromagnetic
energy on $H^{1}(\Omega,\R^2)$:
\begin{align}\label{energy_micro}
E_{\ve}(m) = \ve \int_{\Omega} |D m|^2 + \frac{1}{\ve} \int_{\Omega}(1-|m|^2)^2 + 
\frac{1}{\ve}\int_{\R^2}|\nabla^{-1} \Cu m|^2,
\end{align}
where $m$ is extended by zero outside of $\Omega$. 
They associate to the measurable vector field $m \colon \Omega \to \R^2$ the quantity:
\begin{equation*}
\begin{array}{rccc}
 \chi_m \colon &  \Omega \times \Sp^1 & \longrightarrow     & \left\{0,1 \right\}  \\
             & (x,s)                & \longmapsto & \mathds{1}_{\left\{m(x) \cdot s >0 \right\}}.
\end{array}
\end{equation*}
The introduction of such a quantity is motivated by the following \emph{averaging formula}: let $u \in \Sp^1$,
then
\begin{align}\label{eq:averaging_formula_sphere}
u= \frac{1}{2} \int_{\Sp^1} \mathds{1}_{\left\{u \cdot s>0 \right\}} s \diff \hs^1(s)
\end{align}
An obvious consequence of \eqref{eq:averaging_formula_sphere} is that if we suppose that $|m|=1$ almost everywhere, then for almost every $x$ in $\Omega$,
\begin{align*}
 m(x)= \frac{1}{2} \int_{\Sp^1} \chi_m(x,s) s \diff \hs^1(s).
\end{align*}
They proved then the following theorem:
\begin{thm}[Jabin, Otto, Perthame] \label{thm:J_O_P_energy}
Let $m$ be a zero-state energy, that is a limit of a
sequence $(m_{\ve})_{\ve>0} \subset H^1(\Omega,\R^2)$ with $E_{\ve}(m_{\ve}) \xrightarrow[\ve \to 0]{}  0$. 
Then $|m|=1$ almost everywhere and $m$ satisfies the \emph{kinetic equation}:
\begin{align}\label{eq:kinetic_equation_sphere}
\forall s \in \Sp^1,  \quad & s^{\perp} \cdot  \nabla \chi_m (\cdot,s) = 0 \;
\text{ in } \mathcal{D}'(\Omega).
\end{align}
\end{thm}
Classical results on the regularity of the average of solutions of kinetic equations tell us 
(see \cite{Golse_Lions_Perthame} for example) that equation 
\eqref{eq:kinetic_equation_sphere}
implies for  $m$ to be in $H^{\frac{1}{2}}(\Omega,\R^2)$. But in the case of this particular 
kinetic equation, it turns out that
combining \eqref{eq:kinetic_equation_sphere} with $|m|=1$ implies a much stronger regularity. 
More precisely, the following holds:
\begin{thm}[Jabin, Otto, Perthame]\label{thm:J_O_P_regularity}
Let $m$ be a measurable vector field on $\Omega$ satisfying $|m|=1$ almost everywhere and equation
\eqref{eq:kinetic_equation_sphere}.
Then $m$ is locally Lipschitz continuous inside $\Omega$ except at a locally finite set of singular points. 
Furthermore, 
for every singular point $p$, there is  $\alpha \in \{-1,1\}$ such that in any convex neighborhood
of $p$,
\[
 m(x)=\alpha \frac{x-p}{|x-p|} \quad \text{a.e.}
\]
\end{thm}
Before we state our own results, let us give a bit more insight into the two previous theorems. Looking
at Theorem \ref{thm:J_O_P_energy}, we expect a
zero-state energy $m$ to be of modulus one and curl-free because of the energy 
\eqref{energy_micro}. As a consequence,
the curl-free information about $m$ has to be encompassed into equation \eqref{eq:kinetic_equation_sphere}.
It is indeed easy to see that if a modulus one vector field checks equation \eqref{eq:kinetic_equation_sphere}
then it is curl-free; we will prove this in our more general setting in Proposition
\ref{prop:kinetic_imply_curl}. 
Keeping this in mind, we can know look at Theorem \ref{thm:J_O_P_regularity} in the following way:
 equation \eqref{eq:kinetic_equation_sphere} is a 
\emph{selection principle} among curl-free vector fields taking values into $\Sp^1$ which, 
roughly speaking, does not allow for line-like singularities. 

The question we are addressing here is the following: 
can we find a  kinetic equation which will act as a
selection principle for curl-free vector field taking value into the unit sphere 
of any norm instead of $\Sp^1$? The purpose of this article is to prove that if we make  
simple geometric assumptions on our norm, the answer to this question is yes.
In the spirit of Theorem \ref{thm:J_O_P_regularity}, we will then prove that this
 kinetic equation implies some regularity on the vector field itself. 

Our paper will be divided into three parts. The first one will be devoted to prerequisites on general norms and their modulus of convexity, the second one to the study of a generalized avering formula in the spirit of \eqref{eq:averaging_formula_sphere} and the last one to our main result, i.e. the regularity of solutions to the kinetic equation we are going to introduce.

Before stating our results, we need some definitions on norms and on the way to measure their convexity.

\begin{df}\label{df:convexity_modulus}
Let $\norm{\cdot}$ be a norm and $B$ its closed unit ball. 
We call \emph{modulus of convexity} associated to $B$ the quantity 
$\omega_B \colon \intcc{0}{2} \to \intcc{0}{1}$ defined by:
\begin{align}
\omega_B(\delta) = \inf \left\{ 1-\norm*{\frac{x+y}{2}} : \norm{x}=\norm{y}=1 \text{ and } \norm{x-y} 
\geq \delta \right\}.
\end{align}
\end{df}
\begin{df}\label{df:power_type}
Let $p\geq 0$. We say that a norm is \emph{of power type $p$} if there is a constant $K>0$ such that:
\begin{equation}\label{eq:power_type}
\forall \delta \in [0,2], \qquad  K \delta^p \leq \omega_B(\delta)
\end{equation}
\end{df}
As an example, one should keep in mind that the euclidean norm is of power type $2$ and that more generally,
for $p>1$, $l^p$ norms are of power type $\max \{p,2\}$: this is a consequence of
Clarkson's inequalities and can be found in \cite{Lindenstrauss_Tzafiri_2}. The main reason for introducing theses two notions is that we \emph{need} some property on the norm in order to obtain a kinetic equation which acts as a selection principle excluding line-like singularities.
It turns out that the \emph{power type p} property is exactly the good property to look at.
We refer the reader to the beginning of our third section for a more detailed discussion.

The first important result of our article is an avering formula in the spirit of \eqref{eq:averaging_formula_sphere}:
\begin{thm}\label{thm:representation_formula_convex}
Let $B \subset \R^2$ be a symmetric convex body, that is the unit ball of some norm, then:
\begin{equation}\label{eq:representation_formula_convex}
\forall x \in \partial B, \qquad  x = \frac{1}{2} \int_{\partial B^\perp}
\mathds{1}_{\left\{x \cdot s >0 \right\}} n_{B^{\perp}}(s) \diff \hs^1(s),
\end{equation}
where $n_{B^{\perp}}(s)$ stands for the unit normal to $\partial B^{\perp}$ at $s$.
\end{thm}
\begin{rmk}
 Note that formula \eqref{eq:differentielle_charac} generalizes easily in any dimension in the case
of the euclidean sphere. We refer the reader to \cite{Bochard_Ignat}
for such a formula and its use to prove a rigidity result obtained in the study of an analog
of equation \eqref{eq:kinetic_equation_sphere} in dimension greater than $2$. 
It would be nice to know if we can obtain
an averaging formula like \eqref{eq:representation_formula_convex} in greater dimension for a general 
symmetric convex body.
\end{rmk}

The following corollary of Theorem \ref{thm:representation_formula_convex}, which is crucial to prove our main result, is immediate. 
\begin{cor}\label{cor:representation_formula_field}
Let $\norm{.}$ be a norm, $B$ its closed unit ball and $m \colon \Omega \to \R^2$ a Borel vector field
such that $\norm{m}=1$. We associate to $m$ the quantity:
\begin{equation*}
\begin{array}{rccc}
 \chi_m \colon &  \Omega \times \partial B^{\perp} & \longrightarrow     & \left\{0,1 \right\}  \\
             & (x,s)                & \longmapsto & \mathds{1}_{\left\{m(x) \cdot s >0 \right\}}.
\end{array}
\end{equation*}
Then, for all $x$ in $\Omega$,
\begin{equation}\label{eq:representation_formula_field}
m(x) = \frac{1}{2} \int_{\partial B^\perp} \chi_m (x,s) n_{B^{\perp}}(s) \diff \hs^1(s).
\end{equation}
\end{cor}

We need a last definition to state our main result:

\begin{df}\label{df:vortex}
Let $\norm{.}$ be a norm with unit ball $B$ and $\norm{\cdot}_{\star}$ its dual norm defined as usual by the equality:
\[
 \norm{x}_{\star}=\sup_{\norm{y}\leq 1}\langle x,y \rangle.
\]
We call \emph{vortex} associated to the norm $\norm{\cdot}$ the function $V_B : x \mapsto \nabla_x \norm{\cdot}_{\star}$ at the points where this quantity makes sense. The map $x \mapsto \|.\|$ being convex, $V_B$ is well-defined almost everywhere on $\R^2$ according to Rademacher's theorem.
\end{df}

We can now state the main theorem of our article.
\begin{thm}\label{thm:kinetic_regularization}
Let $\norm{.}$ be a norm of power type $p$  and $B$ its closed unit ball such that $\partial B$ is a
$C^1$ submanifold in $\R^2$. Let  $m \colon \Omega \to \R^2$ be a Borel vector field such
that $\norm{m}=1$ and satisfying:
\begin{align}
\forall s \in \partial B^{\perp},  \quad & n_{B^{\perp}}(s)^{\perp} \cdot  \nabla \chi_m (\cdot,s) 
= 0 \;
\text{ in } \mathcal{D}'(\Omega). \label{eq:kinetic_equation} 
\end{align}
Then $m$ is locally $\frac{1}{p-1}$-H\"older continuous outside a locally finite set of singular point. Furthermore, 
for every singular point $p$, there is  $\alpha \in \{-1,1\}$ such that in any convex neighborhood
of $p$,
\[
 m(x)=\alpha V_B(x-p) \quad \text{almost everywhere},
\]
where $V_B$ is the vortex associated to $\|.\|$ defined in Definition \ref{df:vortex}.
\end{thm}
\begin{rmk}
Notice that the two hypotheses on the norm are not redundant: we can find norms which are not of power type $p$ for any $p$ and whose unit sphere is a $C^1$ submanifold. Take for example  the norm associated to the $l^1$ unit ball whose angles have been smoothed out. On the other hand, an example of a type $p$ norm whose unit ball is \emph{not} a $C^1$ submanifold is given by $\norm{x}=\norm{x}_1+\norm{x}_p$ with $p>2$. The reader interested in the links between convexity and smoothness of the ball of a norm will consult with profit Section 1.e. of \cite{Lindenstrauss_Tzafiri_2}.
\end{rmk}

In order to see how Theorem \ref{thm:kinetic_regularization} could be used, let us note that
Theorem \ref{thm:J_O_P_regularity} itself has interesting consequences. For example, together with clever 
commutator estimates, De Lellis and Ignat proved in \cite{De_Lellis_Ignat} the following:
\begin{thm}[De Lellis, Ignat]\label{thm:regularity_Ignat_DeLellis}
Let $\Omega$ be an open set in $\R^2$ and $r \in [1,3]$. Then every $m$ in
\begin{equation*}
 W_{curl}^{\frac{1}{r},r}(\Omega,\Sp^1):= \left\{m \in  W^{\frac{1}{r},r}(\Omega,\R^2) :
 \mathrm{curl}\, m =0 \text{ in } \mathcal{D}'(\Omega) \text{ and } |m|=1 \text{ a.e.}\right\}
\end{equation*}

is locally Lipschitz continuous inside $\Omega$ except at a locally finite set of singular points,
where $W^{\frac{1}{r},r}(\Omega,\R^2)$ stands for the usual fractional Sobolev space.
\end{thm}
A first natural question would be to ask if we can obtain such a theorem when $\Sp^1$ is changed into the
unit sphere of a norm satisfying the hypotheses of Theorem \ref{thm:kinetic_regularization}.
Another natural question, in view of Theorem \ref{thm:J_O_P_energy}, would be to ask if we can obtain
an analog of Theorem \ref{thm:J_O_P_energy} if we change the energy 
\eqref{energy_micro} into:
\begin{align*}
 \tilde{E}_{\ve}(m) = \ve \int_{\Omega} |D m|^2 + \frac{1}{\ve} \int_{\Omega}(1-\norm{m}^2)^2 + 
\frac{1}{\ve}\int_{\R^2}|\nabla^{-1} \Cu m|^2,
\end{align*}
where $\norm{.}$ is now a norm satisfying the hypotheses of Theorem \ref{thm:kinetic_regularization}.
Note that the proof of both Theorem \ref{thm:regularity_Ignat_DeLellis} and Theorem \ref{thm:J_O_P_energy} 
rely on the notion of \emph{entropy}, that is, in this context, smooth functions $\Phi \colon \Sp^1 \to \R^2$
such that for all open $\Omega$ in $\R^2$,
\[
 \forall m \in C_{curl}^{\infty}(\Omega,\Sp^1), \quad \nabla \cdot [\Phi \circ m]=0.
\]
We refer the reader to \cite{DeLellis_Otto} and \cite{Ignat_Merlet} for more details and 
applications of this notion. Therefore a first step in order to obtain analogs of Theorems 
\ref{thm:J_O_P_energy} and \ref{thm:regularity_Ignat_DeLellis} would be to develop an adapted notion of entropy. 
We will adress this question in a paper to come. 

\section{Modulus of convexity and vortices}\label{sec:convexity}
Our main sources for this section are \cite{Lindenstrauss_Tzafiri_1} and \cite{Lindenstrauss_Tzafiri_2} on 
the geometry of norms and \cite{Ekeland_Temam} on convex analysis. For the reader's comfort, let us collect some general facts on convex functions and convex bodies:
\begin{itemize}
\item A function $f : \R^d \to \R \cup \{+\infty\}$ is said \emph{convex} if $f((1-t)x+ty) \leq (1-t)f(x)+tf(y)$ for all $x,y \in \R^d$, $t \in [0,1]$. Its domain is defined as $\Dom(f) := \{ x : f(x) < \infty\}$.
\item The subdifferential of $f$ at $x$ is defined as $\partial f(x) := \{ p : \forall y, f(y) \geq f(x) + \langle p,y-x\rangle\}$. It is nonempty if $x \in \overset{\circ}{\Dom(f)}$. In that case, $f$ is differentiable at $x$ if and only if its subdifferential is a singleton, in which case $\partial f(x) = \{\nabla f(x)\}$.
\item The Legendre-Fenchel transform of $f$ is defined as
\[f^\star(p) = \sup_x x \cdot p - f(x).\]
If $f$ is a proper ($\not\equiv +\infty$) convex lsc function, so is $f^\star$, and $f = f^{\star\star}$.
\item In that case, $\partial f^\star$ is the inverse of $\partial f$ in the sense that
\begin{equation}\label{eq:subdifferential_inverse}
p \in \partial f(x) \iff x \in \partial f^\star(p).
\end{equation}
\item If $B$ is a convex body and $\norm{\cdot}$ its associated norm, then $\chi_B^\star = \norm{\cdot}_\star$.
\item If $x\in \partial B$, we say that $u$ is a \emph{normal vector to $B$ at $x$} if
\[\forall y \in B, \quad u \cdot (y-x) \leq 0.\]
Notice that the set of normal vectors, denoted by $\mathcal{N}_B(x)$, is a convex cone which is precisely 
$\partial \chi_B(x)$.
\item If $x\in \partial B$, the following relation between the subdifferentials of $\norm{\cdot}$ and $\chi_B$ holds:
\begin{equation}\label{eq:link_subdiff}
\partial_x \norm{\cdot} = \partial \chi_B(x) \cap B^\star, \quad \text{or equivalently} \quad \partial_x \norm{\cdot} = \mathcal{N}_B(x) \cap B^\star.
\end{equation}
\item There is a unique unit normal vector at $x \in \partial B$ if and only if the associated norm is 
differentiable at $x$. In general this is true for $\hs^1$-almost every $x \in \partial B$, and 
it is the case for all $x \in \partial B$ for instance if $\partial B$ is a $\mathcal{C}^1$ submanifold.
\end{itemize} 

Before getting to the heart of the matter, let us explain why the norm should satisfy some property for the kinetic equation to act as a selection principle ruling out line-like singularities. For that purpose, let us take another look at the vortex $V_B$ associated to a norm $\norm{\cdot}$ of unit ball $B$ as defined in Definition \ref{df:vortex}.
It is clear that $V_B$ is curl-free in $\mathcal{D}'(\R^2)$. Indeed, taking $\varphi \in C_0^{\infty}(\R^2)$,
\begin{align*}
\int_{\R^2}\partial_1 \varphi V_{B,2} -  \partial_2 \varphi V_{B,1} \diff x & = 
\int_{\R^2} \nabla \varphi^{\perp} \cdot \nabla \norm{x}_\star \diff x \\
 & = \int_{\R^2} \nabla \cdot \left[ \nabla \varphi^{\perp} \right] \norm{x}_\star \diff x=0,
\end{align*}
where we used the fact that $V_B$ is the distributional gradient of $\norm{\cdot}$ by Rademacher's theorem. 
Furthermore, $V_B$ takes its values into $\partial B$. Indeed by \eqref{eq:subdifferential_inverse}, 
for almost every $x$ (those where $\norm{\cdot}_\star$ is differentiable, \ie where $V_B$ is defined), 
one has $p = V_B(x) \iff x\in \partial_p \norm{\cdot}$, so that $p$ is a point where the subdifferential 
$\partial_p \norm{\cdot}$ contains a nonzero vector, \ie a point on the boundary of the ball. 
Consequently $V_B(x) = p \in \partial B$. We are going to prove that $V_B$ is $\textit{always}$ a 
solution to the kinetic equation \eqref{eq:kinetic_equation}. 
As a consequence, if we want this kinetic equation to rule out line-like singularities then $V_B$ should
\emph{not} exhibit such singularities, and the regularity of $V_B$ is related to the convexity of 
$B$, as we will see in Proposition \ref{prop:convexity_regularity_body}. This is the fundamental reason why 
we need some convexity property on our norm in order to obtain Theorem \ref{thm:kinetic_regularization}.

\subsection{Properties of the normal field}

We start off by stating some useful properties of the normal field 
\begin{equation*}
\begin{array}{rccc}
n_B \colon  \partial B & \longrightarrow &  \mathcal{P}(\Sp^1) \\
       x          & \longmapsto  & \mathcal{N}_B(x) \cap \Sp^1,
\end{array}
\end{equation*}
where $\mathcal{N}_B(x)$ is the normal cone to $B$ at $x$.

\begin{lem}
Let $B$ be the unit ball of a norm and $n_B$ the application defined above. Then, for $\mathcal{H}^1$-almost every $y$ in $\Sp^1$, there is a unique $x$ in $\partial B$ such that $y \in n_B(x)$.
\end{lem}
\begin{proof}
Let us consider: 
 \begin{align*}
&\begin{array}{rccc}
\tilde{n}_B \colon  \partial B & \longrightarrow &  \mathcal{P}(\partial B^{\star}) \\
       x          & \longmapsto  & \mathcal{N}_B(x) \cap \partial B^{\star},
\end{array}&
&\begin{array}{rccc}
\tilde{n}_{B^\star} \colon  \partial B^\star & \longrightarrow &  \mathcal{P}(\partial B) \\
       x          & \longmapsto  & \mathcal{N}_{B^\star}(x) \cap \partial B.
\end{array}
\end{align*}
By \eqref{eq:subdifferential_inverse} and \eqref{eq:link_subdiff}, it is easy to see that $p \in \tilde{n}_B(x)$ if and only if $x \in \tilde{n}_{B^{\star}}(p)$. Since $\partial B^{\star}$ admits a normal vector for $\mathcal{H}^1$-almost every $p$, we get that there is a unique $x$ satisfying $p \in \tilde{n}_B(x)$ for $\mathcal{H}^1$-almost every $p$ in $\partial B^{\star}$. We get uniqueness for $n_B$ by using the fact that the projection from $\Sp^1$ to $\partial B^{\star}$ is a bi-Lipschitz bijection.
\end{proof}
\begin{rmk}
In the rest of the paper, if $y \in \Sp^1$ is such that there is
a unique $x \in \partial B$ satisfying $y \in n_B(x)$, we will abusively 
denote it by $x=n_B^{-1}(y)$. One can deduce from the previous proof that this in in particular the case if $\partial B^\star$ is a $C^1$ submanifold of $\R^2$.
\end{rmk}
\begin{prop}\label{prop:forme_vortex}
Let $\norm{\cdot}$ be a norm of unit ball $B$, then for almost every $x$, 
 \[
  V_B(x) = n^{-1}_B \left(\frac{x}{\abs*{x}}\right).
 \]
\end{prop}

\begin{proof}
We know that $\norm{\cdot}_\star$ is the Legendre-Fenchel transform of $\chi_{B}$. 
Consequently $\norm{\cdot}_\star$ is differentiable at $x$ if and only 
if $x$ is in the normal cone of a unique $p \in \partial B$, that is if and only if 
$n^{-1}_B \left(\frac{x}{\abs*{x}}\right)$ is well defined, in which case 
$n^{-1}_B \left(\frac{x}{\abs*{x}}\right) = p = \nabla_x \norm{\cdot}_\star \doteq V_B(x)$. By the previous lemma, it is true for almost every $x$.
\end{proof}

We can now prove as promised that $V_B$ is always solution to the kinetic equation \eqref{eq:kinetic_equation}.
Indeed, let $\Omega \subset \R^2$, $\varphi \in C_0^{\infty}(\Omega)$, and 
$s \in \partial B$. Using Proposition \ref{prop:forme_vortex}, for almost every $x \in \Omega$,
$V_B(x) = n^{-1}_B \left(\frac{x}{\abs*{x}}\right)$, and
\begin{align*}
\int_{\Omega}\nabla \varphi \cdot n_{B^{\perp}}(s)^{\perp} \chi_{V_B} (x,s) \diff x
=\int_{\Omega \cap \{n^{-1}_B \left(\frac{x}{\abs*{x}}\right)\cdot s>0 \}}\nabla \varphi 
\cdot n_{B^{\perp}}(s)^{\perp} \diff x.
\end{align*}
But for almost every $s$, $\left\{n^{-1}_B \left(\frac{x}{\abs*{x}}\right)\cdot s>0 \right\}$ is 
an half-space whose boundary is the line $L$ passing through $0$ and directed by $n_B(s^{\perp})$, 
whose unit normal is $\pm n_B(s^{\perp})^{\perp}=\mp n_{B^{\perp}}(s)$. 
As a consequence, using Stokes formula and the fact that $\varphi$ is supported in $\Omega$,
\begin{align*}
\int_{\Omega \cap \{n^{-1}_B \left(\frac{x}{\abs*{x}}\right)\cdot s>0 \}}\nabla \varphi 
\cdot n_{B^{\perp}}(s)^{\perp} \diff x = 
\pm n_{B^{\perp}}(s)^{\perp} \cdot  n_{B^{\perp}}(s)\int_{\Omega \cap L} \varphi \diff x=0. 
\end{align*}
Let us define a last quantity that we will need later.
For $u,v \in \R^2$ which are not colinear, we define the convex cone $C(u,v)$ generated by $u,v$ as
\[
 C(u,v):=\left\{ \lambda_1 u+ \lambda_2 v : (\lambda_1,\lambda_2) \in \R^{+} \times \R^{+}  \right\}.
\]

\begin{lem}\label{lem:ordering}
Let $\norm{\cdot}$ be a strictly convex and differentiable (outside $0$) norm of closed unit ball $B$. Then $n_B : \partial B \to \Sp^1$ is monotonic in the sense that
\[w \in C(u,v) \iff n_B(w) \in C(n_B(u),n_B(v)).\]
\end{lem}
\begin{proof}
Notice that $C(u,v)$ is (by Hahn-Banach Theorem) the intersection of the half-spaces which are delimited by vector hyperplanes and which contain $u$ and $v$. First, let us show that, given any such half-space $H$, $n_B$ is a bijection between $\partial B \cap H$ and a half-circle of $\Sp^1$. Since $n_B$ is continuous and $\partial B \cap H$ is compact and connected, $n_B(\partial B \cap H)$ is a compact and connected subset of $\Sp^1$, hence it is a closed arc of $\Sp^1$. Since $\partial B \cap H$ contains two symmetric points $\pm x \in \partial B$, so does $n_B(\partial B \cap H)$, hence it contains a half-circle $[\theta, \theta+\pi]$. It may not be larger than that, otherwise it would contain a nontrivial arc and its symmetric, hence $\partial B \cap H$ as well since $n_B$ is homeomorphic and antisymmetric. This cannot be true. This proves that $n_B$ sends half-spaces of $\partial B$ to half-circles and it is clear that it realizes a bijection between the set of half-spaces of $\partial B$ and the set of half-circles of $\Sp^1$.

Now $w \in C(u,v)$ means that $w$ belongs to any half-space (delimited by a vector hyperplane) containing $u,v$, which means that $n_B(w)$ belongs to any half-circle containing $n_B(u), n_B(v)$ by what we just proved. This exactly means that $n_B(w) \in C(n_B(u), n_B(v))$.
\end{proof}
\subsection{Modulus of convexity}

In this section we give an alternate definition of the modulus of convexity, which offers some advantages: an intuitive interpretation is given in Proposition \ref{prop:power_type_charac}, and a characterization in terms of the regularity of the vortex is given in Proposition \ref{prop:convexity_regularity_body}. This will be useful to prove the sharpness of our main theorem. Of course, it is essentially equivalent to Definition \ref{df:convexity_modulus} which is given in the introduction, as showed in Remark \ref{rk:modulus_equivalence}.

\begin{df}
We define the \emph{modulus of convexity} of $B$ as the greatest nondecreasing function $\rho_B : [0,2] \to \R_+$ such that for all $x,y \in \partial B$ and all $u \in \mathcal{N}_B(x)$,
\begin{equation}\label{eq:convexity_modulus_body}
 u \cdot y \leq u\cdot x - \norm{u}_\star \rho_B(\norm{y-x}).
\end{equation}
We say that $B$ is \emph{of power type $p\geq 2$} if $\rho_B(\delta) \geq C \delta^p$ for some $C > 0$ and that it is \emph{elliptic} if it is of power type 2.
\end{df}
\begin{rmk}
Notice that
\begin{align*}
\forall y \in \partial B,\quad & u \cdot y \leq u\cdot x - \norm{u}_\star 
\rho_B(\norm{y-x}) \\ 
&\text{if and only if} \\
\forall y \in \R^d,\quad & \chi_B(y) \geq \chi_B(x) + u \cdot (y-x) + 
\norm{u}_\star \rho_B(\norm{y-x}).
\end{align*}
This last inequality is reminiscent of the inequality $f(y) \geq f(x) + p \cdot (y-x)$ for $p \in \partial f(x)$ which holds for arbitrary convex functions. The extra term $\norm{u}_\star \rho_B(\norm{y-x})$ measures how much $B$ is convex around its boundary.
\end{rmk}

\begin{rmk}\label{rk:modulus_equivalence}
The two definitions of modulus of convexity are equivalent in the sense that:
\[\rho_B(\delta/2) \leq \omega_B(\delta) \leq \rho_B(\delta)/2.\]
\end{rmk}
\begin{proof}
The proof is quite easy. Take $x,y\in \partial B$ and $u \in \mathcal{N}_B((x+y)/2)$ such that $\norm{u}_\star = 1$. Notice that $u \in \partial_{\frac{x+y}{2}} \norm{\cdot}$ (because $\norm{u}_\star = 1$), hence by convexity:
\[1 = \norm{x} \geq \norm*{\frac{x+y}{2}} + u \cdot \frac{x-y}{2}.\]
By definition of $\rho_B$ one has:
\begin{gather*}
u \cdot y \leq u \cdot \frac{x+y}{2} - \rho_B\left(\norm*{y - \frac{x+y}{2}}\right),
\shortintertext{or equivalently}
u \cdot \frac{y-x}{2} \leq - \rho_B\left(\frac{\norm{y-x}}{2}\right).
 \end{gather*}
Putting this in the previous inequality yields:
\[\norm*{\frac{x+y}{2}} \leq 1 - u \cdot \frac{x-y}{2} \leq 1 - \rho_B\left(\frac{\norm{y-x}}{2}\right).\]
Since $\omega_B$ is precisely the largest map for which the inequality $\norm*{\frac{x+y}{2}} \leq 1 - \omega_B(\norm{y-x})$ holds for all $x,y \in \partial B$, it follows that:
\[ \rho_B (\delta/2) \leq \omega_B(\delta)\]
for all $\delta \in [0,2]$.

Now let us prove the other inequality. This time we take $x,y \in \partial B$, $u \in \mathcal{N}(x)$ with $\norm{u}_\star = 1$. By definition we know that:
\[\norm*{\frac{x+y}{2}} \leq 1 - \omega_B(\norm{y-x}),\]
but since $\norm{u}_\star = 1$, one has:
\[\norm*{\frac{x+y}{2}} \geq u \cdot \frac{x+y}{2} = u \cdot x + u \cdot \frac{y-x}{2} = 1 + u \cdot  \frac{y-x}{2}.\]
Putting these two together, one gets:
\begin{gather*}
u \cdot y \leq u \cdot x - 2\omega_B(\norm{y-x}),
\shortintertext{and by homogenity if $u$ has arbitrary norm:}
u \cdot y \leq u \cdot x - \norm{u}_\star \cdot 2\omega_B(\norm{y-x}).
\end{gather*}
Since $\rho_B$ is the largest function for which such an inequality holds for all $x,y$ and $u \in \mathcal{N}_B(x)$, one has
\[2\omega_B(\delta) \leq \rho_B(\delta)\]
which concludes the proof.
\end{proof}

\begin{figure}[!ht]
\centering
\begin{tikzpicture}[line cap=round,line join=round,>=latex,x=1.0cm,y=1.0cm,scale=0.8]
\clip(-2.5,-2.5) rectangle (2.5,2.5);
\draw (-2,0)-- (0,2);
\draw (0,2)-- (2,0);
\draw (2,0)-- (0,-2);
\draw (0,-2)-- (-2,0);
\draw (0,2)-- (0,-2);
\draw (-2,0)-- (2,0);
\draw [->] (0.28,0.28) -- (0.85,0.85);
\draw [->] (0.28,-0.28) -- (0.85,-0.85);
\draw [->] (-0.28,-0.28) -- (-0.85,-0.85);
\draw [->] (-0.28,0.28) -- (-0.85,0.85);
\begin{scriptsize}
\end{scriptsize}
\end{tikzpicture}
\caption{Unit ball of $\|.\|_1$ and vortex associated to $\|.\|_{\infty}$, singular along the axes.}
\label{pic:unit_ball_l1}
\end{figure}
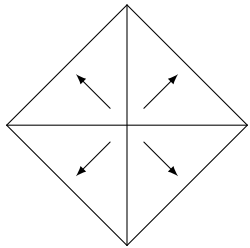

As indicated by its name, the greater $\rho_B$ is, the more convex $B$ is. For example, taking a look at Figure \ref{pic:unit_ball_l1}, it is easy to see that the modulus of convexity associated to the $l^1$ norm is constant equal to $0$ on a neighborhood of $0$. On the other side, the euclidean ball is the most convex of all in the sense that if $\omega_2$ is its modulus of convexity, for all symmetric convex body $B$,
\begin{equation}\label{eq:euclidean_modulus}
\forall \delta \in [0,2], \quad \omega_B(\delta) \leq \omega_2(\delta).
\end{equation}
We refer the reader to \cite{Borne_Norme_2_Nordlander} for a proof of this inequality.
\begin{rmk}
Using the symmetry of $\mathbb{S}^1$, $\rho_2$ can easily be computed and one finds:
\[
 \rho_2(\delta)=1-\sqrt{1-\delta^2} \underset{\delta \to 0}{\sim} \frac{\delta^2}{2}.
\]
As a consequence of this equality, of inequality \eqref{eq:euclidean_modulus} and Remark \ref{rk:modulus_equivalence}, there is \emph{no} norm of power type $p$ with $p<2$.
\end{rmk}

Throughout the rest of this paper, we will use the modulus of convexity $\rho_B$ as defined in this section. The following proposition gives an intuitive characterization of convex bodies of power type $p$. It states that they are convex bodies whose boundary is locally the graph of a map which is above that of the map $t \mapsto \abs{t}^p$ around $0$.

\begin{prop}\label{prop:power_type_charac}
Let $x \in \partial B$, $u$ a normal vector of unit euclidean norm. Take $\gamma : ]-\ve, +\ve[ \to \partial B$ a local parameterization of $\partial B$ around $x = \gamma(0)$, with direct orientation and unit speed. We denote by $\tau = \dot{\gamma}(0)$ the unit tangent at $x$ and $\nu = \tau^\perp$ the inner normal. In the basis $(\tau,\nu)$, we write $\gamma(t) = x + a(t) \tau + b(t) \nu$. Then $B$ is of type $p$ if and only if
\begin{equation}\label{eq:convex_modulus_graph}
b(t) \geq C \abs{a(t)}^p.
\end{equation}
If $B$ is of class $C^2$ and $\kappa(x)$ denotes the curvature at $x$, $B$ is elliptic if and only if for some $C >0$:
\[\forall x \in \partial B, \quad \kappa(x) \geq C.\]
\end{prop}
\begin{proof}
Writing inequality \eqref{eq:convexity_modulus_body} with $y = \gamma(t)$ one gets
\[-b(t) \doteq u \cdot (\gamma(t)-x) \leq -C\abs{\gamma(t)-x}^p \leq -C \abs{a(t)}^p\]
which yields \eqref{eq:convex_modulus_graph}. Conversely, since $b = o(a)$ the inequality $b(t) \geq C \abs{a(t)}^p$ implies \eqref{eq:convexity_modulus_body} for $x$ close to $y$, which is enough to get it for all $x,y$.

If $\partial B$ is of class $C^2$ and so is $\gamma$, by definition one has
\[\gamma(t) = x + t \dot{\gamma}(0) + t^2/2 \cdot \ddot{\gamma}(0) + o(t^2) = x + t \tau + \kappa(x) t^2/2 \cdot \nu + o(t^2).\]
Consequently $a(t) = t + o(t^2)$ and $b(t) = \kappa t^2/2 + o(t^2)$, thus $B$ being elliptic, one has:
\[C(t + o(t^2))^2 \leq \kappa t^2/2 + o(t^2).\]
Dividing by $t^2$ and $t\to 0$ yields $\kappa \geq 2C$ and $\partial B$ has positive curvature. The converse is straightforward.
\end{proof}

The next proposition establishes a link between the power type $p$ property of a norm and the smoothness of its vortex.

\begin{prop}\label{prop:convexity_regularity_body}
Let $B$ be a convex body associated to the norm $\norm{\cdot}$. Then $B$ is of power type $p$ if and only if $\norm{\cdot}_\star$ is differentiable everywhere outside $0$ and $V_B \doteq \nabla \norm{\cdot}_\star$ is $\frac{1}{p-1}$-Hölder continuous far from $0$, in the sense that for some constant $H > 0$ one has:
\[\norm{V_B(v)-V_B(u)} \leq \frac{H}{\delta^{\frac{1}{p-1}}} \norm{u-v}_\star^{\frac{1}{p-1}},\]
for all $\norm{u},\norm{v} \geq \delta$ and all $\delta >0$.
\end{prop}
\begin{proof}
Take $x,y \in \partial B$ and $u\in \mathcal{N}_B(x), v\in \mathcal{N}_B(y)$ with $u,v$ nonzero. One has
\[u\cdot y \leq u \cdot x - \norm{u}_\star\rho_B(\norm{y-x}),\]
\[v\cdot x \leq v \cdot y- \norm{v}_\star\rho_B(\norm{y-x}),\]
and by summing these one gets:
\[2 \min(\norm{u}_\star,\norm{v}_\star) \rho_B(\norm{y-x}) \leq (v-u) \cdot (y-x).\]
Now if $B$ is of power type $p$ then
\[C \norm{y-x}^p \leq \rho_B(\norm{y-x}) \leq \frac{(v-u) \cdot (y-x)}{\min(\norm{u}_\star,\norm{v}_\star)} \leq \frac{\norm{v-u}_\star\norm{y-x}}{\min(\norm{u}_\star,\norm{v}_\star)},\]
thus
\[ \norm{y-x}  \leq \left(\frac{\norm{u-v}_\star}{C\min(\norm{u}_\star,\norm{v}_\star)}\right)^{\frac{1}{p-1}},\]
which yields
\[ \norm{y-x}  \leq \frac{H}{\delta^{\frac{1}{p-1}}} \norm{u-v}_\star^{\frac{1}{p-1}}.\]
provided that $\norm{u}_\star,\norm{v}_\star \geq \delta$, for some constant $H$ depending only on the norm.

Conversely, take $x \in \partial B$, $u \in \mathcal{N}_B(x)$ with $\norm{u}_\star = 1$. We want to show that for all $y\in B$:
\[u\cdot y \leq u \cdot x - C\abs{y-x}^p.\]
Take $x_0$ a minimizer of $g(y) = \chi_B(y) + (u \cdot (x-y) - C\frac{\abs{y-x}^p}{p})$. We want to show that for small enough $C$, $g(x_0) \geq 0$. The necessary optimality condition reads
\[u_0 \coloneqq u + C (x_0-x)^{p-1} \in \partial \chi_B(x_0).\]
Necessarily, $x_0 \in B$ thus
\[\norm{u_0-u}_\star \leq \delta C \norm{x_0-x}^{p-1} \leq \delta C 2^{p-1}\]
where $\delta \coloneqq \sup_{\norm{v} = 1} \norm{v}_\star$.
If $C$ is small enough, say $C \delta 2^{p-1} \leq 1/2$ then both $\norm{u}_\star, \norm{u_0}_\star \geq 1/2$. On the other hand one has by hypothesis:
\[\norm{x_0-x}^{p-1} \leq 2 H^{p-1} \norm{u_0-u}_\star \leq 2 \delta C H^{p-1}\norm{x_0-x}^{p-1}.\]
Now if one takes $C$ even smaller, so that $2\delta C H^{p-1} < 1$, one must have $x=x_0$ and $u = u_0$, which implies that $\min_{\R^2} g = g(x_0) = g(x) = 0$ and the following holds:
\[\forall y \in B,\quad u\cdot y \leq u \cdot x - C\abs{y-x}^p/p \qquad \text{with } C = \min(2^{-p}, 1/(2\delta H^{p-1})).\]
If $u$ has arbitrary norm, this becomes:
\[\forall y \in B, \quad  u\cdot y \leq u \cdot x - C' \norm{u}_\star \norm{y-x}^p\]
where $C' = C/p$.
\end{proof}

As a corollary, one gets the following regularity theorem.

\begin{thm}\label{thm:regularity_of_normal}
Let $\norm{\cdot}$ be a norm with closed unit ball $B$. It is of power type $p$ if and only if $n_{B}^{-1} : \Sp^1 \to \partial B$ is well-defined and $\frac{1}{p-1}$-Hölder continuous. In particular it is elliptic iff $n_{B}^{-1}$ is Lipschitz continuous.
\end{thm}

\section{Averaging formula}\label{sec:averaging_formula}

This section is devoted to the proof of Theorem \ref{thm:representation_formula_convex}. A first interesting 
remark is that if we want an averaging formula in the spirit of formula \eqref{eq:representation_formula_convex} 
to be true, we \emph{need} the convex to be symmetric. 
More precisely, we have the following proposition:
\begin{prop}\label{prop:representation_formula_gene}
Let $B \subset \R^2$ be a convex body with $0$ in its interior. Assume there is a Borel vector measure 
$\mu \in \cM^2(\R^2\setminus \{0\})$ such that the following
representation formula holds:
\begin{align}\label{eq:representation_symmetry}
\forall x \in \partial B, \qquad x =  \int_{\R^2} \mathds{1}_{\{ x \cdot s >0 \}} \mu(\diff s).
\end{align}
Then $B$ is symmetric and one may replace $\mu$ in the formula by the measure $\bar{\mu}$ defined by
\[\bar{\mu}(A) \coloneqq \frac{\mu(A) - \mu(-A)}{2}\]
for all Borel set $A$, which is an antisymmetric measure in the sense that $\bar{\mu}(-A) = -\bar{\mu}(A)$, and which 
does not give mass to vector lines.
\end{prop}
\begin{proof}
Let us set
\begin{align}\label{eq:x+_def}
x^+ &\coloneqq \{ s \in \R^2 : \langle x, s \rangle > 0\},& 
x^- &\coloneqq \{ s \in \R^2 : \langle x, s \rangle < 0\},&
x^{\perp} &\coloneqq \{s \in \R^2  : \langle x,s \rangle =0 \}.
\end{align}
One may define the gauge of $B$ as
\begin{equation*}
\begin{array}{lcll}
j_B \colon  & \R^2 & \longrightarrow & \R_+ \\
            & x   & \longmapsto & \inf \{ t> 0 : x \in t B\},
\end{array}
\end{equation*}
the radial projection onto $\partial B$ and radial symmetry with respect to $B$ respectively as:
\begin{align*}
\begin{array}{lcll}
p_B \colon &  \R^2 \setminus \{0\} & \longrightarrow & \partial B \\
           &  x   & \longmapsto & \frac{x}{j_B(x)},
\end{array}\qquad
\begin{array}{lcll}
s_B \colon & \R^2 \setminus \{0\} & \longrightarrow & \R^2 \\
           &  x   & \longmapsto & j_B(x) p_B(-x).
\end{array}
\end{align*}
First, let us assume that $\mu$ does not charge lines of the plane, \ie $\mu(D) = 0$ 
for all vector line $D \subset \R^2$.
Writing $\R^2 = x^+ \sqcup x^- \sqcup x^\perp$ and noticing that $s_B(x)^+ = x^-$, we obtain:
\[
\mu(x^+)+\mu(s_B(x)^+)=\mu(x^+) + \mu(x^-)+\mu(x^{\perp})=  \mu(\R^2) \eqqcolon v \in \R^2 .
\]
If $x$ is in $\partial B$, according to \eqref{eq:representation_symmetry} one has 
$x = \mu(x^+)$ and $s_B(x)=\mu(s_B(x)^+)$ since $s_B(x) \in \partial B$, which yields
\begin{equation}\label{eq:sym_relation}
\forall x \in \partial B, \qquad x+s_B(x)=v.
\end{equation}
But $s_B(x)$ is colinear to $x$, and since $\overset{\circ}{B} \neq \emptyset$, we can find two non
colinear vector $x_1, x_2 \in \partial B$.
Now writing \eqref{eq:sym_relation} with $x_1$ and $x_2$ implies $v=0$; consequently, 
for all $x \in \partial B$, $x=-s_B(x)=-p_B(-x)$, which exactly means that $B$ is symmetric.

Now we may get rid of the hypothesis that $\mu$ does not charge lines. If $\mu$ gives a positive mass to some vector lines, it may only charge countably many of them since $\mu$ is finite and these lines only intersect at $0$ which is $\mu$-negligible. Let us denote by $L$ the reunion of the perpendiculars to these lines. The previous reasoning shows that $x + s_B(x) = v$ for all $x \in \partial B \setminus L$. But $\partial B \cap L$ is countable hence $\partial B \setminus L$ is dense in $\partial B$ and by continuity of $s_B$ the identity
\[x + s_B(x) = 0 = v\]
holds for all $x \in \partial B$, which implies that $B$ is symmetric and $\mu$ does not charge vector lines. Thus if $x \in \partial B$ then $-x \in \partial B$, which implies:
\begin{gather*}
-x = \int_{\R^2} \mathds{1}_{\{-x \cdot s >0\}} \mu(\diff s) = \mu(-(x^+)) \shortintertext{and}
x = \frac{x - (-x)}{2} = \frac{\mu(x^+) - \mu(-(x^+))}{2} := \bar{\mu}(x^+).
\end{gather*}
Obviously, $\bar{\mu}$ is antisymmetric and does not charge lines since $\mu$ does not either.
\end{proof}

Now let us pass on to the proof of the representation formula as stated in Theorem \ref{thm:representation_formula_convex}. It turns out that checking such a formula is quite easy once one has found the right candidate for $\mu$, but the following theorem is more precise as it states that $\mu$ is essentially unique.

\begin{thm}\label{thm:representation_formula_convex_complete}
Let $B \subset \R^2$ be a symmetric convex body, that is the unit ball of some norm, 
then there is a unique antisymmetric measure $\mu$ supported on $\partial B^\perp$ such that:
\begin{equation}
\forall x \in \partial B, \qquad  x = \int_{\R^2}
\mathds{1}_{\left\{x \cdot s >0 \right\}} \mu(\diff s),
\end{equation}
which is
\[\mu = 1/2 \cdot n_{B^\perp} \hs^1_{\mres \partial B^\perp}\]
where $n_{B^{\perp}}(s)$ stands for the unit normal to $\partial B^{\perp}$ at $s$.
\end{thm}

\begin{proof}
Let $B$ be a symmetric convex body and suppose that there is $\mu \in \mathcal{M}^2(\R^2)$ supported in 
$\partial B^\perp$ and satisfying the representation formula:
\begin{align}
\forall x \in \partial B, \qquad x =  \int_{\partial B^\perp} \mathds{1}_{\{x \cdot s>0\}} \mu(\diff s).
\end{align}
We want to prove that $\mu = \frac{n_{B^\perp}}{2}  \hs^1_{\mres \partial B^\perp}$.
 Recalling the notations of \eqref{eq:x+_def}, notice that 
this rewrites more concisely as
\begin{equation}\label{eq:rep_concise}
\forall x \in \partial B, \qquad x = \mu(x^+).
\end{equation}
If $u,v$ are two points of $\partial B^\perp$, we denote by $\intoc{u}{v}$ the oriented arc of $\partial B^\perp$ 
delimited by $u$ and $v$ ($v$ included and $u$ not included).  Let us set $x = u^\perp,y = v^\perp$, which both 
belong 
to $\partial B$ since it is symmetric.
\begin{figure}\label{pic:mon_cul_sur_la_commode}
\centering
 \begin{tikzpicture}[line cap=round,line join=round,x=1.0cm,y=1.0cm]
\clip(-2.5,-2.5) rectangle (2.5,2.5);
\draw [rotate around={25.14:(0,0)}] (0,0) ellipse (1.5cm and 1.12cm);
\draw [rotate around={-64.86:(0,0)},line width=0.4pt] (0,0) ellipse (1.5cm and 1.12cm);
\draw (0,0)-- (-0.08,1.42);
\draw (0,0)-- (-0.64,1.36);
\draw (0,0)-- (-1.42,-0.08);
\draw (-1.36,-0.64)-- (0,0);
\begin{scriptsize}
\fill [color=black] (0,0) circle (0.6pt);
\draw[color=black] (0.16,-0.1) node {$0$};
\draw[color=black] (-0.75,1.5) node {$v$};
\draw[color=black] (0,1.6) node {$u$};
\draw[color=black] (-1.6,0) node {$x$};
\draw[color=black] (-1.5,-0.8) node {$y$};
\draw[color=black] (1.8,0.5) node {$\partial B$};
\draw[color=black] (-1.3,1.3) node {$\partial B^{\perp}$};
\end{scriptsize}
\end{tikzpicture}
\caption{}
\end{figure}

Notice that if the oriented angle between $u$ and $v$ is such that 
$0 < \sphericalangle(u,v) < \pi$ then $]u,v] = x^+ \setminus y^+$. Finally we set $E := \intoc{u}{v}$ and 
$F := x^+ \setminus E$, so that
\begin{align*}
x^+ &= E \sqcup F,&
y^+ &= F \sqcup (-E),\shortintertext{thus}
\mu(x^+) &= \mu(E) + \mu(F),&
\mu(y^+) &= \mu(F) + \mu(-E).
\end{align*}
Since $\mu$ is antisymmetric, $\mu(-E) = -\mu(E)$ hence substracting the last equality to the previous one yields 
$\mu(x^+) - \mu(y^+) = 2\mu(E)$. 
Using \eqref{eq:rep_concise} one gets
\[x - y = 2 \mu(\intoc{u}{v}), \qquad \ie \qquad \mu(\intoc{u}{v}) = - \frac{1}{2} (v^\perp - u^\perp) =
\frac{1}{2}R^{-1} (v-u).\]
In short, setting $F = \frac{1}{2} R^{-1}$, for all $u,v$ such that $0 < \sphericalangle(u,v) < \pi$, we have:
\begin{equation}\label{eq:stielt}
\mu(\intoc{u}{v}) = F(v) - F(u).
\end{equation}
This relation\footnote{Notice that \eqref{eq:stielt} would match the definition of the Stieltjes measure 
associated to a given function 
$F$ if $\mu$ was a measure on the real line. In this case, $\mu$ is given by the distributional derivative of 
$F$ and we are 
merely transposing this fact to our case, where $\mu$ is a measure on a closed curve.} allows us to assert that 
$\mu$ is the following 
measure on $\partial B^\perp$:
\begin{equation}\label{eq:stielt_charac}
\mu(\diff s) = D_s F \tau(s) \diff s = F  \tau(s) \diff s
\end{equation}
where $\tau(s)$ is the unit vector tangent to $\partial B^\perp$ at $s$ (with direct orientation). 
By definition, $n_{B^\perp}(s) = R^{-1} \tau(s)$ hence \eqref{eq:stielt_charac} rewrites 
$\mu = \frac{n_{B^\perp}}{2} \hs^1_{\mres \partial B^\perp}$, which is what we shall 
prove now. Setting $\nu = n_{B^\perp} \hs^1_{\mres \partial B^\perp}$, it is easy to check that 
$\nu(\intoc{u}{v}) = F(v) - F(u) = \mu(\intoc{u}{v})$ for all $u,v$ such that 
$0 < \sphericalangle(u,v) < \pi$. Indeed, take a curve $\gamma : \intcc{0}{l} \to L^\perp$ which 
describes the oriented arc $\intcc{u}{v}$ and is parameterized by arc-length. Then one has
\begin{align*}\label{eq:rep_charac_bis}
\begin{split}
\nu(]u,v]) = \int_0^l n_{B^\perp}(\gamma(t)) \frac{\diff t}{2} &= 
\int_0^l R^{-1}(\dot{\gamma}(t)) \frac{\diff t}{2}\\
&= \frac{1}{2} R^{-1}(\gamma(l) - \gamma(0))\\
&= F(v) - F(u).
\end{split}
\end{align*}
Since the Borel sets $\intoc{u}{v}$ form a $\pi$-system generating the Borel $\sigma$-algebra of 
$\partial B^\perp$, $\mu = \nu$. Moreover, since $\nu$ satisfies relation \eqref{eq:stielt}, 
it satisfies the representation formula \eqref{eq:rep_concise} by taking $v \to -u$. This is what we 
wanted to prove.
\end{proof}

\section{Kinetic formulation}\label{holder_regularity}
\subsection{Kinetic formulation and curl-free vector fields of norm one}
This subsection is devoted to justifying the heuristic explanation of our introduction, 
that is that our kinetic equation acts as a selection principle for norm one curl-free vector fields.
First, we have the following easy proposition:
\begin{prop}\label{prop:kinetic_imply_curl}
Let $\Omega$ be an open set in $\R^2$, $B$ the closed unit ball associated to the norm $\norm{\cdot}$ and $m$ a Borel vector field such that $\norm{m}=1$. If $m$ satisfies the kinetic equation \eqref{eq:kinetic_equation},
then $m$ is curl-free in the distributional sense.
\end{prop}
\begin{proof}
Let $\varphi \in C_0^{\infty}(\Omega)$. Then,
\begin{align*}
\int_{\Omega} \partial_1 \varphi m_2 -\partial_2 \varphi m_1 \diff x  = 
\int_{\Omega} \nabla \varphi^{\perp} \cdot m \diff x & 
=   \int_{\Omega} \nabla \varphi^{\perp} \cdot \left(\frac{1}{2}\int_{\partial B^\perp} \chi_m(x,s)
n_{B^\perp}(s) 
\diff \hs^1(s) \right) \diff x  \\
& = \frac{1}{2} \int_{\partial B^\perp} \int_{\Omega} \nabla \varphi \cdot  n_{B^\perp}(s)^{\perp}  
\chi_m(x,s) \diff x \diff \hs^1(s)=0,
\end{align*}
where we use the Corollary \ref{cor:representation_formula_field} in the first line and the kinetic equation
\eqref{eq:kinetic_equation} in the second, which concludes the proof.
\end{proof}
The next proposition gives kind of a reciprocal if the field $m$ is smooth.
\begin{prop}\label{prop:curl_smooth_imply_kinetic}
Let $\Omega$ be an open convex set in $\R^2$, 
 $x \mapsto \norm{x}$ a norm in $C^2(\R^2 \setminus \{0\})$, and 
 $m \in C^1(\Omega,\R^2)$ a curl-free vector 
field such that $\norm{m}=1$. Then $m$ satisfies the kinetic equation \eqref{eq:kinetic_equation}.
\end{prop}
\begin{proof}
Let $m$ be a vector field satisfying our hypothesis. Using the characteristics method, it is easy to see 
that such a field is constant along lines. Indeed, writing $F(x):=\norm{x}$, differentiating the relation 
$F(m(x))=1$ in $\Omega$ and using the fact that $m$ is curl-free, one gets:
\begin{equation}\label{eq:differentielle_charac}
D_x m \cdot \nabla F(m(x))= D_x^T m \cdot \nabla F(m(x))=0.
\end{equation}
Now let $t \mapsto \gamma(t)$ be a smooth curve taking values in $\Omega$ and $p(t):=m(\gamma(t))$.  
Differentiating $p$, one obtains $\dot{p}(t)=D_{\gamma(t)} m \cdot \dot{\gamma}(t)$. Therefore, 
choosing $\gamma$
to be a solution of the differential equation $\dot{\gamma}(t)= \nabla F(m(\gamma(t)))$ and 
using \eqref{eq:differentielle_charac} the couple $(\gamma, p)$ is now solution of the system:
\begin{align*}
\left\{
\begin{array}{l}
 \dot{p}(t)=0, \\
 \dot{\gamma}(t)= \nabla F(p(t)).
\end{array}
\right.
\end{align*}
Thus $p$ is constant equal to $m(\gamma(0)) = m(x)$, $\gamma$ satisfies $\dot{\gamma} = \nabla F(m(x))$ and $m$ is constant along $\gamma$, that is $m$ is constant on $L_x \cap \Omega$ where $L_x$ is the line passing through $x$ and directed by $\nabla F(m(x))$. Set $s \in \partial B^\perp$. 
If $m(x_0) \cdot s = 0$, then $m(x_0) = \pm s^\perp$, say $m(x_0) = s^\perp$. 
Since $m$ is constant along the line $\gamma : t \mapsto x_0 + t\nabla F(s^\perp)$ and equal to $s^\perp$, then $\chi_m(\gamma(t),s) = 0$  and by differentiation at $t = 0$ one gets $\nabla F(s^\perp) 
\cdot \nabla \chi_m(x_0,s) = 0$ (we write $\xi \cdot \nabla f$ to designate the partial derivative
of $f$ in the direction $\xi$ provided it exists). But $\nabla F(s^\perp)$ is the normal vector 
$\xi$ to $B$ at $s^\perp$ such that $\norm{\xi}_\star = 1$, 
in particular it is colinear to $n_{B}(s^\perp)$, thus 
$n_{B}(s^\perp)\cdot \nabla \chi_m(x_0,s) = 0$. 
Now if $m(x_0) \cdot s > 0$, then $m(x) \cdot s > 0$ along $\gamma$, 
otherwise it must vanish somewhere along that line and $m \equiv \pm s^\perp$ on that line, which 
is a contradiction. Consequently, $\chi_m(\cdot,s)$ is constant along $\gamma$, thus we get again by 
differentiation along $\gamma$ at $t=0$:
\[
n_{B}(s^\perp)\cdot \nabla \chi_m(x_0,s) = 0.
\]
The same reasoning works for $m(x_0) \cdot s < 0$. Noticing that 
$n_{B^\perp}(s)^\perp = n_{B}(s^\perp)$, this implies that
\begin{align}\tag{\ref{eq:kinetic_equation}}
\forall s \in \partial B^{\perp}, \forall x \in \Omega, \quad & n_{B^\perp}(s)^{\perp} \cdot  
\nabla \chi_m (x,s) = 0.
\end{align}
\end{proof}

\subsection{Direction conservation and trace theorem}
We are now entering the heart of the proof of Theorem \ref{thm:kinetic_regularization}, which is quite similar 
in spirit to the proof of Theorem \ref{thm:J_O_P_regularity} by Jabin Otto and Perthame. 
In order to use equation \eqref{eq:kinetic_equation} to gain regularity, an important point 
will be to be able to define the trace of $m$ along a line. 
\begin{rmk}
Note that starting from here, we will \emph{always} place ourselves under the asumptions of Theorem 
\ref{thm:kinetic_regularization}, that is the unit sphere $\partial B$ associated to a norm will be
a $C^1$ submanifold and the norm itself will be of power type $p$.
\end{rmk}

\begin{df}
Let $m  \in \mathcal{L}^1_{\mathrm{loc}}(\Omega,\R^2)$. We say that $x_0$ is a \emph{Lebesgue point} of $m$ if:
\[
 \lim_{r \to 0} \fint_{B_r(x_0)}|m(x)-m(x_0)| \diff x=0.
\]
We denote by $\Leb(m) \subset \Omega$ the set of Lebesgue point of $m$. It is well known 
(see \cite{Evans_Gariepy} for example) that $|\Omega \setminus \Leb(m)|=0$.
\end{df}

Let us state a useful lemma which relates Lebesgue points of $m$ to those of $\chi_m(\cdot,s)$.
\begin{lem}\label{lem:lebesgue_points}
Let $\Omega \subset \R^2$ be an open set and $m \in \mathcal{L}^1_{\mathrm{loc}}(\Omega,\partial B)$. 
\begin{enumerate}[label=(\roman*)]
\item If $x_0$ is a Lebesgue point of $\chi_m(\cdot,s)$ for $\hs^1$-almost every $s \in \partial B^\perp$ 
then it is a Lebesgue point of $m$.\\ \label{lebesgue_points_item1}
\item If $x_0$ is a Lebesgue point of $m$ then $x_0$ is a Lebesgue point of $\chi_m(\cdot,s)$ for all
$s \in \partial B$ such that $m(x_0) \cdot s \neq 0$ .\label{lebesgue_points_item2}
\end{enumerate}
In particulier $x_0$ is a Lebesgue point of $m$ if and only if it is a Lebesgue point of $\chi_m(\cdot,s)$ for almost every $s$.
\end{lem}
\begin{proof}
The proof of \ref{lebesgue_points_item1} is a simple consequence of the averaging formula
\eqref{eq:representation_formula_field}.
Indeed, if $x_0$ is a Lebesgue point of $\chi_m(\cdot,s)$ for $\mathcal{H}^1$-a.e. $s \in \partial B^\perp$:
\begin{align*}
\fint_{B_r(x_0)} \abs*{m(x) - m(x_0)} \diff x  &= \frac{1}{2}\fint_{B_r(x_0)} \abs*{\int_{\partial B^\perp} (\chi_m(x,s) - \chi_m(x_0,s)) n_{B^\perp}(s) \diff \mathcal{H}^1(s)} \diff x\\
& \leq \frac{1}{2} \fint_{B_r(x_0)}\int_{\partial B^\perp} \abs{\chi_m(x,s)-\chi_m(x_0,s)} \diff \mathcal{H}^1(s) \diff x \\
& = \frac{1}{2} \int_{\partial B^\perp}  \fint_{B_r(x_0)} \abs{\chi_m(x,s)-\chi_m(x_0,s)} \diff x  \diff \mathcal{H}^1(s) \xrightarrow{r\to 0} 0,
\end{align*}
by the dominated convergence theorem.

Now let us prove \ref{lebesgue_points_item2}. Take a Lebesgue point $x_0$ of $m$ and 
$s \in \partial B^\perp$ such that 
$ m(x_0) \cdot s  \neq 0$, say for example $ m(x_0) \cdot s > 0$. First, let us remark that
\begin{align*}
 m(x_0) \cdot s  \frac{\abs{B_r(x_0) \cap \{  m \cdot s  \leq 0\}}}{\abs{B_r(x_0)}}
& \leq \frac{1}{B_r(x_0)} \int_{B_r(x_0) \cap \{ m \cdot s  \leq 0\}}   (m(x_0) \cdot  s -   m(x) \cdot  s) 
\diff x\\
& \leq \fint_{B_r(x_0)} \abs{m(x_0) - m(x)} \diff x.
\end{align*}
Dividing by $ m(x_0) \cdot s $ and taking the limit when $r$ goes to $0$, one gets
\[
\frac{\abs{B_r(x_0) \cap \left\{  m \cdot s \leq 0  \right\}}}{\abs{B_r(x_0)}} \xrightarrow{r \to 0} 0,
\]
but
\begin{align*}
\fint_{B_r(x_0)} \abs{\chi_m(x,s) - \chi_m(x_0,s)} \diff x &= \fint_{B_r(x_0)} \abs{\chi_m(x,s) - 1} \diff x\\
&=  \frac{\abs{B_r(x_0) \cap \{  m \cdot s  \leq 0\}}}{\abs{B_r(x_0)}} \xrightarrow{r \to 0} 0,
\end{align*}
hence $x_0$ is a Lebesgue point of $\chi_m(\cdot,s)$. The case $ m(x_0) \cdot s<0$ is done in a similar way.
\end{proof}

We can now prove the following proposition, which plays a crucial role in the proof of Theorem 
\ref{thm:kinetic_regularization}.
\begin{prop}\label{prop:Propagation_scalar_product}
Let $\Omega \subset \R^2$ be a convex open set and $m \in \mathcal{L}^1_{loc}(\Omega,\partial B)$ satisfying 
the kinetic equation \eqref{eq:kinetic_equation}. Assume that $y,z \in \Leb(m)$ 
and that $s \in \partial B^{\perp}$ with $z-y \in \vect( n_{B^\perp}(s)^{\perp})$. Then
\begin{align}
m(y)\cdot s >0 \quad (resp. <0) & \implies m(z) \cdot s \geq 0 \quad (resp. \leq 0). \label{scalar_product_pos} 
\end{align}
\end{prop}
\begin{proof} 
Let us suppose that $m(y) \cdot s \neq0$. If $m(z) \cdot s =0$, then the proof is over; assume otherwise that 
$m(z)\cdot s \neq 0$. Fix $x \in \Omega$ and define for $\ve$ small enough:
\begin{gather*}
\chi^\ve(x) \coloneqq \fint_{B_\ve(x)} \chi_m(y,s) \diff y,
\shortintertext{which may be written as}
\chi^\ve(x) = \rho^\ve * \chi_m(\cdot,s) \quad \text{where} \quad \rho^\ve \coloneqq 
\frac{\indic_{B_{\ve}(0)}}{\abs{B_{\ve}(0)}}.
\end{gather*}
$\chi^\ve$ is continuous and for every Lebesgue point $x_0$ of $\chi_m(\cdot, s)$, 
one has $\chi^\ve(x_0) \xrightarrow{\ve \to 0} \chi_m(x_0,s)$.
Now let $(\delta_n)_{n\in \N^*}$ be a sequence of smooth mollifiers supported in $B_{1/n}$ and set
\[
 \chi_n^\ve \coloneqq \delta_n * \chi^\ve.
\]
We are going to show that $\chi_n^\ve(y)=\chi_n^\ve(z)$ and for that we introduce
\begin{align*}
 \begin{array}{rccc}
g \colon &  [0,1]  & \longrightarrow     & \R  \\
         &     t   & \longmapsto & \chi_n^\ve(y + t(z-y)).
\end{array}
\end{align*}
The function $g$ is smooth and satisfies
\begin{align*}
g'(t) &= (z-y) \cdot \nabla \chi_n^\ve(y+t(z-y))\\
&= (z-y) \cdot \brc[\Big]{\nabla \sqbrc{\delta_n * \rho^\ve} * \chi_m(\cdot,\xi)}\brc*{y+t(y-z)}\\
&= \int_{\R^2} (z-y) \cdot \nabla \left[\delta_n * \rho_{\ve} \right] (y+t(z-y)-v) \chi_m(v,s) \diff v.
\end{align*}
But $(z-y) \in \text{Vect}( n_{B^\perp}(s)^{\perp})$ and $\delta_n * \rho_{\ve} \in 
\cC_c^{\infty}(\Omega)$; 
therefore because of \eqref{eq:kinetic_equation}, $g'(t)=0$ and
\[
\chi_n^\ve(y)=g(0)=g(1)=\chi_n^\ve(z).
\]
Taking the limit as $n$ goes to $+\infty$, we obtain
\[
\chi^\ve(y)=\chi^\ve(z).
\]
But $y,z$ are Lebesgue points of $m$ and $m(y) \cdot s \neq 0, m(z) \cdot s \neq 0$
hence by \ref{lebesgue_points_item2} of Lemma \ref{lem:lebesgue_points} $y,z$ are
Lebesgue points of $\chi_m(\cdot,s)$, thus taking the limit as $\ve$ goes to $0$, we have 
\[
\chi_m(y,s)=\lim_{\ve \to 0} \chi^\ve(y) = \lim_{\ve \to 0} \chi^\ve(z)=\chi_m(z,s),
\]
which concludes the proof.
\end{proof}

\begin{thm}\label{thm:trace}
Let $\Omega \subset \R^2$ be a convex open set and $m \in \mathcal{L}_{loc}^1(\Omega,\partial B)$ 
satisfying the kinetic 
equation \eqref{eq:kinetic_equation}. Let $L$ be a segment of 
the form $L = \{ t v : t \in [-1,1]\}$ for some $v \in \Sp^1$. For $r>0$ set
\[
P_r:=\{x_1v^{\perp}+x_2v : (x_1,x_2) \in [-r,r] \times [-1,1] \}
\]
and assume that $P_r$ is included in $\Omega$ for small $r$. Then there exists a measurable function 
$\tilde{m} \colon [-1,1] \to \R^2$ satisfying
\begin{enumerate}[label=(\roman*)]
\item $\displaystyle \lim\limits_{r \to 0}  \fint_{P_r} |m(x)-\tilde{m}(x_2)| \diff x =0$,\label{trace1}
\item for almost every $x_2$, $\tilde{m}(x_2) \in \partial B$,\label{trace2}
\item if $x \in \Leb(m) \cap L$ then $x_2 \in \Leb(\tilde{m})$ and $m(x) = \tilde{m}(x_2)$,\label{trace3}
\item if $x_2 \in \Leb(\tilde{m})$, then there exists a sequence $x_n \in \Leb(m)$ such that
\[x_n \to x \quad \text{and} \quad m(x_n) \to \tilde{m}(x_2).\] \label{trace4}
\end{enumerate}
The vector field $\tilde{m}$ is called the trace of $m$ on the line $L$.
\end{thm}
\begin{rmk}
We only assume the form of $L = \{ t v : t \in [-1,1]\}$ for commodity but the result obviously holds 
for any general segment in $\Omega$. This will be used several times in the rest of the article. 
\end{rmk}

Equation \eqref{eq:kinetic_equation} means that $\chi_m$ only depends on $n_{B^\perp}(s)\cdot x$ 
and $s$. The next lemma asserts that $\chi_m$ may be written as a Borel function of these quantities.
\begin{lem}\label{lem:borel_factor}
If $m$ satisfies \eqref{eq:kinetic_equation}, for any convex set $V$ 
(for instance a convex neighborhood of $L$) there exists a Borel function $\tilde{\chi}_m : 
\R \times \partial B^\perp \to \{0, 1\}$ such that for a.e. $x$ in $V$,
\begin{align}\label{eq:trace_un_d}
 \chi_m(x,s)=\tilde{\chi}_m(n_{B^\perp}(s)\cdot x,s).
\end{align}
\end{lem}

\begin{proof}
We know that $m$ and $\chi_m$ are Borel maps. For $s$ fixed, we denote by 
$x = u_1 n_{B^\perp}(s) + u_2 n_{B^\perp}(s)^\perp$ the decomposition of 
$x$ on $n_{B^\perp}(s), n_{B^\perp}(s)^\perp$. Taking arbitrary smooth functions $\phi,\psi$ 
of the variables $u_1,u_2$ with (suitable) compact support and testing the kinetic equation against 
$\rho(x) = \phi(u_1)\psi(u_2)$, one gets
\[\int \phi(u_1) \psi'(u_2) \chi_m(x,s) \diff u_2 \diff u_1 = 0,\]
which implies that for a.e. $u_1$, the map $u_2 \mapsto 
\chi_m(u_1 n_{B^\perp}(s) + u_2 n_{B^\perp}(s)^\perp,s)$ is a.e. equal to a constant, namely to 
the mean value of $\chi_m(\cdot,s)$ over the slice $V_{s,x\cdot n_{B^\perp}(s)}$ where 
$V_{s,a} =  \{ a n_{B^\perp}(s) + t n_{B^\perp}(s)^\perp : t \in \R\} \cap V$.
Therefore, setting
\[
\tilde{\chi}_m(a,s) = \fint_{V_{s,a}} \chi_m(x,s)\diff x,
\]
one has
\[
\chi_m(x,s) = \tilde{\chi}_m(n_{B^\perp}(s) \cdot x, s)
\]
for all $s$ and almost every $x$ in $V$. Let us justify briefly that $\tilde{\chi}_m$ is Borel. 
Setting $F(s,a,t) = a n_{B^\perp}(s) + t n_{B^\perp}(s)^\perp$, which is continuous, and 
$U_{s,a} = \{t :  F(s,a,t) \in V\}$, one has
\[\tilde{\chi}_m(a,s) = \frac{1}{\abs{U_{s,a}}} 
\int_{U_{s,a}} \chi_m(a n_{B^\perp}(s) + t n_{B^\perp}(s)^\perp,s) \diff t = 
\frac{\int_\R \chi_m(F(s,a,t),s) \mathds{1}_V(F(s,a,t)) \diff t}{\int_\R \mathds{1}_V(F(s,a,t)) \diff t},\]
and it is now clear that it is Borel as an integral of a Borel map with parameter.
\end{proof}

We can now prove Theorem \ref{thm:trace}.

\begin{proof}[Proof of Theorem \ref{thm:trace}]
We proceed in several steps.

\paragraph{\textbf{Step 1:} Trace for $\chi_m$.\\}
Provided that $v \cdot n_{B^\perp}(s) \neq 0$, the function $[-1,1] \ni x_2 \mapsto 
\tilde{\chi}_m(x_2 v \cdot n_{B^\perp}(s),s)$ is the trace of $\chi_m(\cdot,s)$ on $L$ in the sense that:
\begin{equation}\label{eq:trace_chi_und}
\lim_{r \to 0} \fint_{P_r}
\abs*{\chi_m(x,s)- \tilde{\chi}_m(x_2 v \cdot n_{B^\perp}(s),s)} \diff x =0.
\end{equation}
Indeed, using Lemma \ref{lem:borel_factor}:
\begin{align*}
2 \fint_{P_r} 
|\chi_m(x,s)- \tilde{\chi}_m(& n_{ B^\perp}(s)\cdot v x_2,s)| \diff x   \\ 
& =\int_{-1}^1 \fint_{-r}^r
\abs*{\chi_m(x_1 v^{\perp}+x_2 v,s)- \tilde{\chi}_m(x_2 v \cdot n_{ B^\perp}(s)  ,s)} 
\diff x_1 \diff x_2 \\
& =\int_{-1}^1 \fint_{-r}^r
\abs*{\tilde{\chi}_m((x_1 v^{\perp} +x_2 v)\cdot n_{ B^\perp}(s),s)- 
\tilde{\chi}_m(x_2 v \cdot n_{B^\perp}(s),s)} \diff x_1 \diff x_2 \\
& \leq \frac{1}{|v \cdot n_{B^\perp}(s)|} \int_{-1}^1 \fint_{-r}^r
\abs*{\tilde{\chi}_m(x_1 v^{\perp}\cdot n_{B^\perp}(s)+y_2,s)- 
\tilde{\chi}_m(y_2,s)} \diff x_1 \diff y_2 \\
& \leq \frac{1}{|v \cdot n_{B^\perp}(s)|} \sup_{|y_1| \leq r} \int_{-1}^1 
\abs*{\tilde{\chi}_m(y_1+y_2,s)- 
\tilde{\chi}_m(y_2,s)} \diff y_2. \\
\end{align*}
This last quantity goes to $0$ as $r$ goes to $0$ because of the continuity of the translation in $L^1$.
Note that for a fixed $v \in \Sp^1$, there are only two vectors $s \in \partial B^{\perp}$ 
for which $v \cdot n_{B^\perp}(s) = 0$ because $B$ is a convex body which is stricly convex (the norm being of type $p$). As a consequence equality 
\eqref{eq:trace_chi_und} is true for $\mathcal{H}^1$-a.e. $s$ in $\partial B^{\perp}$.
\paragraph{\textbf{Step 2:} Trace for $m$.\\}
For $x_2 \in [-1,1]$, we define the trace $\tilde{m}$ of $m$ on $L$ 
by the following equality :
\[
 \tilde{m}(x_2) =  
 \frac{1}{2} \int_{\partial B^\perp} \tilde{\chi}_m(x_2 v \cdot n_{B^\perp}(s)  ,s) 
 n_{B^\perp}(s) \diff \mathcal{H}^1(s).
\]
Then, thanks to the averaging formula \eqref{eq:representation_formula_field},
\begin{align*}
 2\fint_{P_r} \abs*{m(x)-\tilde{m}(x_2)} \diff x 
 & \leq \fint_{P_r} \int_{\partial B^\perp} 
 \abs*{(\chi_m(x,s)- \tilde{\chi}_m(x_2 v \cdot n_{B^\perp}(s),s))n_{B^\perp}(s)} 
 \diff \mathcal{H}^1(s) 
 \diff x\\
 & \leq  \int_{\partial B^\perp} \left( \fint_{P_r} 
 \abs*{\chi_m(x,s)- \tilde{\chi}_m(x_2 v \cdot n_{B^\perp}(s),s)}  \diff x \right) \diff \mathcal{H}^1(s).
\end{align*}
Using \eqref{eq:trace_chi_und} and the dominated convergence theorem, \ref{trace1} is proved. Moreover,
\begin{align*}
\int_{-1}^1 \left| 1-\norm{\tilde{m}(x_2)} \right| \diff x_2  & \leq 2 \fint_{P_r} 
\abs*{\norm{m(x)}-\norm{\tilde{m}(x_2)}} \diff x \\
& \leq  C \fint_{P_r}|m(x)-\tilde{m}(x_2)| \diff x,
\end{align*}
the last equality resulting from the equivalence of norms on $\R^2$ and the reverse 
triangle inequality. This last quantity going to $0$ as $r$ goes to $0$, \ref{trace2} is proved.

\paragraph{\textbf{Step 3:} Proof of \ref{trace3}.\\}
Take $x$ a Lebesgue point of $m$. We know by Lemma \ref{lem:lebesgue_points} that $x$ is a Lebesgue point of $\chi_m(\cdot,s)$ for almost all $s$. We define the cube 
\[Q_r^s(x):=x + \{a n_{B^\perp(s)} + b n_{B^\perp(s)}^\perp : (a,b) \in [-r,r]^2 \},\]
and write $y = y_1^s n_{B^\perp(s)} + y_2^s n_{B^\perp(s)}^\perp$ for arbitrary $y$. Notice that:
\[\fint_{Q_r^s(x)} \abs{\chi_m(y,s)-\chi_m(x,s)} \diff y  = \fint_{Q_r^s(x)} \abs{\tilde{\chi}_m(y_1^s,s)-\tilde{\chi}_m(x_1^s,s)} \diff y = \fint_{x_1^s-r}^{x_1^s+r} \abs{\tilde{\chi}_m(u,s)-\tilde{\chi}_m(x_1^s,s)} \diff u.\]
The quantity on the left tends to $0$ as $r$ goes to $0$ hence the quantity on the right as well, which means that $x_1^s = x \cdot n_{B^\perp}(s) \in \Leb(\tilde{\chi}_m)$ for almost all $s$. Moreover, by definition of the trace, one has for all $y \in L$:
\[\tilde{m}(y_2) =  
 \frac{1}{2} \int_{\partial B^\perp} \tilde{\chi}_m(y_2 v \cdot n_{B^\perp}(s)  ,s) 
 n_{B^\perp}(s) \diff \mathcal{H}^1(s),\]
thus
\[ \fint_{x_2-r}^{x_2+r} \abs{\tilde{m}(y_2) - \tilde{m}(x_2)} \diff y_2 \leq \int_{\partial B^\perp} \underbrace{\fint_{x_2-r}^{x_2+r} \abs{\tilde{\chi}_m(y_2 v \cdot n_{B^\perp}(s),s)- \tilde{\chi}_m(x_2 v \cdot n_{B^\perp}(s),s)}\diff y_2}_{J_s(r)} \diff s.\]
Provided that $v \cdot n_{B^\perp}(s) \neq 0$, which is true for a.e. $s$:
\[J_s(r) = \fint_{x_1^s - r v \cdot n_{B^\perp}(s)}^{x_1^s + r v \cdot n_{B^\perp}(s)} \abs{\tilde{\chi}_m(u,s)- \tilde{\chi}_m(x_1^s,s)}\diff u \xrightarrow{r \to 0} 0,\]
hence by the dominated convergence theorem, $J_s(r)$ being bounded by $2$, one gets:
\[ \fint_{x_2-r}^{x_2+r} \abs{\tilde{m}(y_2) - \tilde{m}(x_2)} \diff y_2 \to 0,\]
and $x_2$ is a Lebesgue point of $\tilde{m}$. Moreover:
\begin{align*}
\tilde{m}(x_2) & = \lim_{r \to 0} \frac{1}{2}\fint_{x_2-r}^{x_2+r} \tilde{m}(t) \diff  t \\ 
&= \lim_{r \to 0} \frac{1}{2} \int_{\partial B^\perp} \fint_{x_2-r}^{x_2+r}  \tilde{\chi}_m(t v \cdot n_{B^\perp}(s)  ,s) \diff t \;n_{B^\perp}(s)\diff s\\
&= \lim_{r \to 0} \frac{1}{2}\int_{\partial B^\perp} \fint_{x_2-r}^{x_2+r} \fint_{x_1-r}^{x_1+r} \tilde{\chi}_m((t v + u n_{B^\perp}(s)^\perp)\cdot n_{B^\perp}(s),s) \diff u \diff t \;n_{B^\perp}(s)\diff s\\
&= \lim_{r \to 0} \frac{1}{2}\int_{\partial B^\perp} \fint_{x_2-r}^{x_2+r} \fint_{x_1-r}^{x_1+r} \chi_m(t v + u n_{B^\perp}(s)^\perp,s) \diff u \diff t \;n_{B^\perp}(s)\diff s\\
&= \lim_{r \to 0} \frac{1}{2} \int_{\partial B^\perp} \fint_{L_r^s(x)} \chi_m(y,s) \diff y\; n_{B^\perp}(s) \diff s
\end{align*}
where we have defined the parallelogram (which is non-flat for a.e. $s$):
\[L_r^s(x) = x + \{ t v + u n_{B^\perp}(s) : (t,v) \in [-r,r]\}.\]
Since $x$ is a Lebesgue point of $m$, it is a Lebesgue point of $\chi_m(\cdot,s)$ for almost all $s$,
hence 
\[
\fint_{L_r^s(x)} \chi_m(y,s) \diff y \to \chi_m(x,s).
\]
Passing to the limit as $r$ goes to $0$ in the integral:
\[\tilde{m}(x_2) = \frac{1}{2} \int_{\partial B^\perp} \chi_m(x,s) n_{B^\perp}(s) \diff s = m(x),\]
using again the representation formula.

\paragraph{\textbf{Step 4:} Proof of \ref{trace4}.\\}
Fix $\ve > 0$, and take $\delta >0$ to be suitably chosen later in terms of $\ve$. For $a,r$ positive and smaller than $\ve$, we look at the quantity $E(a,r) = \fint_{x_2 -a}^{x_2+a} \fint_{-r}^r \abs{m(y)-\tilde{m}(x_2)} \diff y$. This quantity is bounded from above as follows:

\begin{align*}
E(a,r) = \fint_{x_2 -a}^{x_2+a} \fint_{-r}^r \abs{m(y)-\tilde{m}(x_2)} \diff y \leq &
\fint_{x_2 -a}^{x_2+a} \fint_{-r}^r \abs{m(y)-\tilde{m}(y_2)} \diff y  \\ 
+ \fint_{x_2 -a}^{x_2+a} \fint_{-r}^r \abs{\tilde{m}(y_2)-\tilde{m}(x_2)} \diff y\\
\leq & \fint_{P_{r,a}(x)} \abs{m(y)-\tilde{m}(y_2)} \diff y + \fint_{x_2 -a}^{x_2+a} \abs{\tilde{m}(y_2)-\tilde{m}(x_2)} \diff y_2,
 \end{align*}
where
\[
P_{r,a}:=x+ \{u_1v^{\perp}+u_2v : (u_1,u_2) \in [-r,r] \times [-a,a]\}.
\]
Since $x_2$ is a Lebesgue point of $\tilde{m}$, one may find $a$ small enough such that
\[\fint_{x_2-a}^{x_2+a} \abs{\tilde{m}(y_2) - \tilde{m}(x_2)} \diff y_2 \leq \delta.\]
Then for this fixed $a$, by \ref{trace1} the quantity $\fint_{P_{r,a}} \abs{m(y)-\tilde{m}(y_2)} \diff y$ goes to $0$ as $r$ goes to $0$, thus one may find $r$ small such that it is less then $\delta$, yielding $E(a,r) \leq 2\delta$. Now $E(a,r)$ may be bounded from below as follows:
\[E(a,r) = \fint_{P_{r,a}} \abs{m(y)-\tilde{m}(x_2)} \geq \ve \frac{\abs{y \in \Leb(m) \cap P_{r,a}: \abs{m(y)-\tilde{m}(x_2)}> \ve}}{\abs{P_{r,a}}},\]
which implies that 
\[\frac{\abs{y \in \Leb(m) \cap P_{r,a}: \abs{m(y)-\tilde{m}(x_2)}> \ve}}{\abs{P_{r,a}}} \leq \frac{2\delta}{\ve}.\]
Choosing $\delta = \ve / 4$ this quantity is strictly less than $1$ hence the set
\[\{y \in \Leb(m) \cap P_{r,a}: \abs{m(y)-\tilde{m}(x_2)} \leq \ve\}\]
is not empty. Any point $x^\ve$ in this set is a Lebesgue point such that $\abs{x^\ve -x} \leq \sqrt{2} \ve$ (because $x^\ve \in P_{\ve,\ve}(x)$), and sastisfying $\abs{m(x^\ve)-\tilde{m}(x_2)}\leq \ve$. Taking $\ve = 1/n$ gives the desired conclusion.
\end{proof}

With this definition of the trace, we obtain as a corollary an extension of Proposition \ref{prop:Propagation_scalar_product}.
\begin{cor}\label{cor:prod_scalaire_trace}
Assume that $\Omega \subset \R^2$ is a convex open set and $m \in \mathcal{L}^1_{loc}(\Omega,\partial B)$ satisfies the kinetic equation \eqref{eq:kinetic_equation}. If $L$ is a given line segment in $\Omega$ and $z\in L$ is such that $z_2$ is a Lebesgue point of the trace $\tilde{m}$ on $L$, then for all $y \in \Leb(m)$ such that $z-y \in \vect( n_{B^\perp}(s)^{\perp})$ one has:
\begin{align}
m(y)\cdot s >0 \quad (resp. <0) & \implies \tilde{m}(z_2) \cdot s \geq 0 \quad (resp. \leq 0), \label{scalar_product_pos_trace}.
\end{align}
\end{cor}
\begin{proof}
The key point lies in \ref{trace4} of Theorem \ref{thm:trace}. With the notations of this theorem, if $z$ is a point of $L$ such that $z_2$ is a Lebesgue point of $\tilde{m}$, one may find a sequence $z^n \in \Leb(m)$ such that $z^n \to z$ and $m(z^n) \to \tilde{m}(z_2)$. We know that $(z-y) \cdot n_{B^\perp}(s) = 0$. The map $n_{B^\perp} : \partial B^\perp \to \Sp^1$ being a homeomorphism, one may find a sequence $s^n \in \partial B^\perp$ such that $(z^n-y) \cdot n_{B^\perp}(s^n) = 0$ and $s^n \to s$ (set $s^n = \pm n_{B^\perp}^{-1}\left(\frac{(z^n-y)^\perp}{\abs{z^n-y}}\right)$). Applying Proposition \ref{prop:Propagation_scalar_product}, one gets
\[m(z^n) \cdot s^n \geq 0,\]
and passing to the limit $n\to \infty$ yields $\tilde{m}(z_2) \cdot s \geq 0$.
\end{proof}

\subsection{The regularity theorem}
Now that we have defined a trace of our vector field along lines, our goal is to prove a form of invariance
along lines just like what would happen in the smooth case  thanks to the characteristic method.
\begin{prop}\label{prop:invariance_droite}
Let $m \in \mathcal{L}^1(\Omega,\partial B)$ satisfying the kinetic equation \eqref{eq:kinetic_equation}. Suppose 
that $x_0$ is a Lebesgue point of $m$ and denote by $L$ be the line passing trough $x_0$ and directed by 
$v:= n_{B}(p) \in \Sp^1$ where $m(x_0) = p$. Then for almost every $x_2$, 
$\tilde{m}(x_2) \in \{\pm p\}$ where $\tilde{m}$ is the trace of 
$m$ on $L$ in the sense of Theorem \ref{thm:trace}.
\end{prop}
\begin{proof}
By translation of the domain and rotation of the target space, we may assume without loss of generality 
that $x_0 = 0$ and $p = (\lambda,0)$ with $\lambda > 0$. Let $h>0$, $r_0>0$ be fixed and set
\[
E_{r_0}^{\pm}:= \left\{ x \in \Omega \cap \Leb(m) : x=h v\pm t v^{\perp} \quad 0< t\leq r_0  \right\}.
\]
\begin{figure}[!ht]
\centering
\begin{tikzpicture}[line cap=round,line join=round,x=1.0cm,y=1.0cm]
\clip(-1,-1) rectangle (6,4);
\draw [->,>=latex] (0,0) -- (1.58,0);
\draw [domain=-1.59:6.16] plot(\x,{(-0--2.18*\x)/4.3});
\draw [thick] (4.64,1.51)-- (3.96,2.85);

\draw[<->] (4.23,2.98)-- (4.57,2.32);
\draw[<->] (-0.14,0.27)-- (4.16,2.45);

\draw[dashed] (0,0)-- (-0.14,0.27);
\draw[dashed] (3.96,2.85)-- (4.23,2.98);

\draw (4.21,2.13)-- (4.26,2.04);
\draw (4.35,2.09)-- (4.26,2.04);

\begin{scriptsize}
\draw[color=black] (0,-0.2) node {$O$};
\draw[color=black] (0.86,-0.15) node {$p$};
\draw[color=black] (2.2,1.72) node {$h$};
\draw[color=black] (5.7,2.6) node {$L$};
\draw[color=black] (4.7,2.7) node {$r_0$};
\draw[color=black] (3.7,3) node {$E_{r_0}^+$};
\draw[color=black] (5,1.6) node {$E_{r_0}^-$};
\end{scriptsize}
\end{tikzpicture}
\caption{The set $E^{\pm}_{r_0}$}
\end{figure}
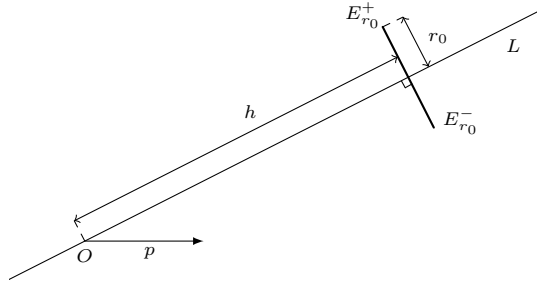

For $x^{\pm} \in E_{r_0}^{\pm}$, we will note  
$s^{\pm}=(s_1^{\pm},s_{2}^{\pm}):=n_{B^{\perp}}^{-1}\left( \frac{(x^{\pm})^{\perp}}{|x^{\pm}|} \right) \in \partial B^{\perp}$
and $\alpha^{\pm}=(\alpha_1^{\pm},\alpha_2^{\pm}):=n_{B^{\perp}}^{-1} \left( \frac{hv^{\perp} \mp r_0 v}{|hv^{\perp} \mp r_0 v|} 
\right)$. Remark that $s^{\pm}$ is chosen in a way such that 
$x^{\pm} \in \text{Vect}(n_{B^{\perp}}(s^{\pm})^{\perp})$.
We claim that for $r_0$ small enough, $s_2^{\pm} \geq \frac{\lambda}{2}>0$ and there is $C>0$ such that:
\begin{align*}
(0>s_1^+ \geq C \alpha_1^+  &\text{ and } 0<s_1^- \leq C \alpha_1^-) \\
& \text{ or } \\ 
(0<s_1^+ \leq C \alpha_1^+  &\text{ and } 0>s_1^- \geq C \alpha_1^-),
\end{align*}
Let us prove this claim. First note that $n_{B^{\perp}}(s^{\pm})$ is in the cone 
$C(v^{\perp}, n_{B^{\perp}}(\alpha^{\pm}))$ as defined above Lemma \ref{lem:ordering}. 
Indeed,
\[
 n_{B^{\perp}}(s^{\pm})=\frac{(x^{\pm})^{\perp}}{|(x^{\pm})^{\perp}|}=\frac{hv^{\perp} \mp t v}{|hv^{\perp} \mp t v|},
\]
and
\begin{align*}
\frac{1}{|hv^{\perp} \mp t v|}\left( hv^{\perp}\mp tv \right)
 & =\frac{1}{|hv^{\perp} \mp t v|} \left( \frac{t}{r_0}(h v^{\perp}\mp r_0v)+ 
 h \left(1- \frac{t}{r_0} \right)  v^{\perp} \right),\\
 & = \frac{1}{|hv^{\perp} \mp t v|} \left( \frac{t|h v^{\perp}\mp r_0v|}{r_0}n_{B^{\perp}}(\alpha^{\pm})+ 
 h \left(1- \frac{t}{r_0} \right)  v^{\perp} \right).
\end{align*}
Using Lemma \ref{lem:ordering},
this implies $s^{\pm} \in C(n_{B^{\perp}}^{-1}(v^{\perp}), 
\alpha^{\pm})$. 
But $n_{B^{\perp}}^{-1}(v^{\perp})=n_{B}^{-1}(v)^{\perp}=p^{\perp}$ 
so that:
\[
s^{\pm} \in C(p^{\perp},\alpha^{\pm}),
\]
meaning that there are $t_1,t_2$ in $\R^{+*}$ such that
\begin{align*}
\begin{pmatrix}
s_1^{\pm} \\
s_2^{\pm}
\end{pmatrix}
=t_1
\begin{pmatrix}
0 \\
\lambda
\end{pmatrix}
+t_2 \alpha^{\pm}.
\end{align*}

Because of the continuity of $n_{B^{\perp}}^{-1}$ proved in Theorem \ref{thm:regularity_of_normal},
$\alpha^{\pm} \xrightarrow[]{r_0 \to 0} p^{\perp}=(0,\lambda)$,
which implies that for $r_0$ small enough $s_2^{\pm} \geq \frac{\lambda}{2}>0$. 
To conlude the proof of our claim, we just have to remark that there are only two possibilities for 
the position of $\alpha^{\pm}$ with respect to $p^{\perp}$ since $p^\perp \in C(\alpha^+,\alpha^-)$ by Lemma \ref{lem:ordering}.
\begin{figure}[!ht]
\begin{tikzpicture}[line cap=round,line join=round,>=latex,x=1.0cm,y=1.0cm]
\clip(-1,-0.5) rectangle (4,2);
\draw [->] (0,0) -- (0.6,0.52);
\draw [->] (0,0) -- (0,1);
\draw [->] (0,0) -- (-0.81,1.02);
\draw [->] (3,0) -- (3,1);
\draw [->] (3,0) -- (3.6,0.52);
\draw [->] (3,0) -- (2.19,1.02);
\begin{scriptsize}
\fill [color=black] (0,0) circle (0.5pt);
\draw[color=black] (0.3,0.8) node {$p^{\perp}$};
\draw[color=black] (-0.4,1) node {$\alpha^+$};
\draw[color=black] (0.8,0.4) node {$\alpha^-$};
\fill [color=black] (3,0) circle (0.5pt);
\draw[color=black] (3.3,0.8) node {$p^{\perp}$};
\draw[color=black] (-0.4+3,1) node {$\alpha^-$};
\draw[color=black] (0.8+3,0.4) node {$\alpha^+$};
\end{scriptsize}
\end{tikzpicture}
\caption{Two differents configurations}
\end{figure}
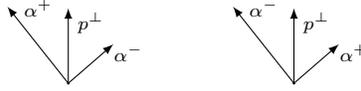

In the first configuration, because $s^{\pm}, \alpha^{\pm} \in \partial B^{\perp}$, 
$s^{+} \in C(p^{\perp},\alpha^{+})$ implies $0>s_1^+ \geq C \alpha_1^+ $ and 
$s^{-} \in C(p^{\perp},\alpha^{-})$ implies $0<s_1^- \leq C \alpha_1^-$. 
In the second configuration, using the same trick ones gets 
$s^{+} \in C(p^{\perp},\alpha^{+})$ implies $0<s_1^+ \leq C \alpha_1^+ $ and 
$s^{-} \in C(p^{\perp},\alpha^{-})$ implies $0>s_1^- \geq C \alpha_1^-$. 
This concludes the proof of our claim.

Without loss of generality, we will place ourselves in the second case that is we will chose $r_0$ small enough
such that $s_2^{\pm} \geq \frac{\lambda}{2}>0, 0<s_1^+ \leq C \alpha_1^+$ and $0>s_1^- \geq C \alpha_1^-$.
Then, using Proposition \ref{prop:Propagation_scalar_product}, for $x^+ \in E_{r_0}^+$,
\begin{align*}
  m(0) \cdot s^+=\lambda s_1^+>0 & \implies m(x^+)\cdot s^+=m_1(x^+) s_1^+ +m_2(x^+)s_2^+ \geq 0 \\
                   & \implies m_2(x^+)\geq -\frac{m_1(x^+) s^+_1}{s^+_2} \geq -\tilde{C}\alpha^+_1.
\end{align*}
Identically, we get for $x^- \in E_{r_0}^-$:
\begin{align*}
  m(0) \cdot s^-=\lambda s_1^-<0 & \implies m(x^-)\cdot s^-=m_1(x^-) s_1^- +m_2(x^-)s_2^- \leq 0 \\
        & \implies m_2(x^-)\leq -\frac{m_1(x^-) s^-_1}{s^-_2} \leq -\tilde{C} \alpha_1^-.
\end{align*}
Notice that Theorem \ref{thm:trace} implies that
\begin{align*}
\int_{-1}^1 \fint_{0}^r \abs*{m(x)-\tilde{m}(x_2)} \diff x_1 \diff x_2 = 2 \fint_{P_r} \abs*{m(x)-\tilde{m}(x_2)} \diff x \xrightarrow{r\to 0} 0
\shortintertext{and}
\int_{-1}^1 \fint_{-r}^0 \abs*{m(x)-\tilde{m}(x_2)} \diff x_1 \diff x_2 = 2 \fint_{P_r} \abs*{m(x)-\tilde{m}(x_2)} \diff x \xrightarrow{r\to 0} 0.
\end{align*}
Consequently there is a sequence $r_n \to 0$ such that:
\[
\fint_{0}^{r_n} \abs*{m(x)-\tilde{m}(x_2)} \diff x_1\xrightarrow{n \to +\infty} 0 \quad \text{and} \quad 
\fint_{-r_n}^0 \abs*{m(x)-\tilde{m}(x_2)} \diff x_1 \xrightarrow{n \to +\infty} 0
\]
for almost every $x_2$. Moreover we know that for almost every $x_2$, $x_1 v^{\perp} +x_2 v$ is a Lebesgue point of $m$ for almost every $x_1$, so that $m_2(x_1 v^{\perp} +x_2 v) \geq -\tilde{C}\alpha_1^+$ if $x_1 \geq 0$ and $m_2(x_1 v^{\perp} +x_2 v) \leq -\tilde{C}\alpha_1^-$ if $x_1 \leq 0$. Thus it follows that for almost every $x_2$:
\begin{align*}
\tilde{m}_2(x_2)= \lim _{n\to\infty} \fint_{0}^{r_n} m_2(x_1 v^{\perp} +x_2 v) \diff x_1 \geq 
-\tilde{C} \alpha^+_1
                   \xrightarrow{r \to 0} 0,
\shortintertext{and}
\tilde{m}_2(x_2)= \lim _{n\to\infty} \fint_{-r_n}^0 m_2(x_1 v^{\perp} +x_2 v) \diff x_1 \leq 
-\tilde{C} \alpha^-_1 
                   \xrightarrow[]{r \to 0} 0. 
\end{align*}
This shows that $\tilde{m}_2(x_2)=0$ for a.e. $x_2$, i.e. $\tilde{m}(x_2) \in \{ \pm p\}$ which is what 
we wanted to prove.
\end{proof}

We can now finish the proof of Theorem \ref{thm:kinetic_regularization}.
\begin{proof}[Proof of Theorem \ref{thm:kinetic_regularization}]
Consider any open convex $\omega \subset \Omega$ with $d=d( \omega, \partial \Omega) > 0$. It is enough to prove Theorem \ref{thm:kinetic_regularization} in any such $\omega$. 
Let $x,y$ in $\omega \cap \Leb(m)$, $L_x$ be the line passing through $x$ and directed by $n_B(m(x))$
and $L_y$ be the line passing through $y$ and directed by $n_B(m(y))$. 
If  $L_x$ and $L_y$ are parallel and distinct, we claim that $m(x)=m(y)$. Indeed, if it is not the case,
$m(x)=-m(y)$. We choose $s \in \partial B^{\perp}$ such that $y-x \in \text{Vect}(n_{B^{\perp}}(s)^{\perp})$.
Note that $m(y)\cdot s \neq 0$, otherwise $s$ and $m(y)^{\perp}$ are colinear and  
$y-x \in \text{Vect}(n_{B}(m(y)))$
and $L_x=L_y$. Now using Proposition \ref{prop:Propagation_scalar_product} 
$\text{sign}(m(y)\cdot s)= \text{sign}(m(x) \cdot s)$ hence $m(x)=m(y)$. 

If $L_x$ and $L_y$ intersect each other, let $O$ be their intersection point. 
Up to a change of variable in the target space,
we can suppose that $O$ is the origin, 
$L_x= \text{Vect}(\frac{x}{|x|})$ and $L_y= \text{Vect}(\frac{y}{|y|})$. Now, up
to a change of sign for $m$, we can suppose that $m(x)=n_B^{-1}(\frac{x}{|x|})$ 
and $m(y)= \pm n_B^{-1}(\frac{y}{|y|})$.
We are going to prove that the trace of $m$ along $L_x$ and $L_y$ in the sense of Theorem 
\ref{thm:trace} is completely 
determined by $m(x)$ and $m(y)$. More precisley,
\begin{align*}
\text{for a.e. $t \in \R \text{ such that } t\frac{x}{|x|} \in L_x \cap \omega$}, \quad \tilde{m}_{L_x}(t)= n^{-1}_B \left( \frac{x}{|x|} \right), \\
\text{for a.e. $t \in \R \text{ such that } t\frac{y}{|y|} \in L_y \cap \omega$}, \quad \tilde{m}_{L_y}(t)= n^{-1}_B \left( \frac{y}{|y|} \right),
\end{align*}
where $\tilde{m}_{L_x}$ and $\tilde{m}_{L_y}$ stand for the trace of $m$ along $L_x$ (resp $L_y$) in the sense of 
Theorem \ref{thm:trace}.
Let $z=t\frac{y}{|y|} \in L_y \cap \omega$ chosen such that $t$ is a Lebesgue point of $\tilde{m}_{L_y}$
and choose $s \in \partial B^{\perp}$ such that $z-x \in \text{Vect}(n_{B^{\perp}}(s)^{\perp})$. Because of
Proposition \ref{prop:invariance_droite}, 
$\tilde{m}(t)=\pm m(y)= \pm n_{ B}^{-1}\left( \frac{y}{|y|} \right)$.
But
\begin{align*}
 m(x) \cdot s & = n_{ B}^{-1} \left( \frac{x}{|x|} \right) 
 \cdot n_{ B^{\perp}}^{-1}\left( \frac{(z-x)^{\perp}}{|z-x|} \right) 
  = n_{B}^{-1} \left( \frac{x}{|x|} \right) 
 \cdot n_{ B}^{-1}\left( \frac{(z-x)}{|z-x|} \right)^{\perp}.
\end{align*}
Because $z$ and $x$ are not colinear, $\frac{(z-x)}{|z-x|} \not\in \text{Vect} (x)$ and 
$n_{B}^{-1}$ being bijective, $m(x) \cdot s \neq 0$. Because of Proposition \ref{cor:prod_scalaire_trace},
this implies $\text{sign}(\tilde{m}_{L_y}(t)\cdot s)= \text{sign}(m(x) \cdot s)$. But $\frac{y}{|y|}=\frac{z}{|z|} \in C \left(\frac{x}{|x|}, 
\frac{z-x}{|z-x|} \right)$ and
by Lemma \ref{lem:ordering}, this implies 
$n^{-1}_B \left(\frac{y}{|y|} \right) \in 
C\left(n^{-1}_B\left(\frac{x}{|x|} \right),n_{B}^{-1} \left( \frac{z-x}{|z-x|} \right) \right)$, that is there is 
$\alpha,\beta$ in $\R^+$ such that:
\[
 n^{-1}_B \left(\frac{y}{|y|} \right)=\alpha n^{-1}_B\left(\frac{x}{|x|} \right) + 
 \beta n_{B}^{-1} \left( \frac{z-x}{|z-x|} \right).
\]
Now,
\begin{align*}
 n^{-1}_B \left(\frac{y}{|y|} \right)\cdot s & =n^{-1}_B \left(\frac{y}{|y|} \right) \cdot n_{B}^{-1} \left( \frac{z-x}{|z-x|} \right)^{\perp} \\
 & = 
 \alpha n^{-1}_B\left(\frac{x}{|x|} \right)  \cdot n_{B}^{-1} \left( \frac{z-x}{|z-x|} \right)^{\perp} = 
 \alpha m(x) \cdot s,
\end{align*} 
so that $\text{sign}(\tilde{m}_{L_y}(t)\cdot s)= \text{sign}(n^{-1}_B \left(\frac{y}{|y|}\right) \cdot s)$. 
This forces $\tilde{m}_{L_y}(t)=n^{-1}_B \left(\frac{y}{|y|}\right)$ and we proved that $\tilde{m}_{L_y}$ is 
determined along $L_y$. We show in the same way that $\tilde{m}_{L_x}=n^{-1}_B \left(\frac{x}{|x|}\right)$ along $L_x$. 
Note that if we started with $m(x)=-n_B^{-1}\left(\frac{x}{|x|} \right)$ instead of 
$n_B^{-1} \left(\frac{x}{|x|} \right)$, the trace of 
$m$ along $L_x$ and $L_y$ would be determined in the same way, but with the opposite sign. We will use this
in the end of our proof.

We are now going to distinguish two cases: \\
\textbf{Case 1}: $\text{dist}(O,\omega)\leq d$. We claim that $m$ is a vortex in $\omega$.
First note that up to a change of sign for $m$, we can suppose that $m(x)=n_B^{-1} \left(\frac{x}{|x|} \right)$
and $m(y)=n_B^{-1} \left(\frac{y}{|y|} \right)$.
Let us take $z \in \omega \cap \Leb(m) \cap (L_x \cup L_y)^c$ and let $L_z$ be the line directed 
by $m(z)$ and passing through $z$. 
If we prove that $L_z$ intersects $L_x$ and $L_y$ in $O$, we are done. Suppose now that this 
is not the case. Using the same argument
than in the beginning, $L_z$ cannot be parallel neither to $L_x$ nor to $L_y$.
Let $P_1=L_z \cap L_x$ and $P_2=L_z \cap L_y$. By convexity, there is a small portion of 
$L_z$ in $\Omega$. But using the remark in our first part to the lines $L_x$ and $L_z$, one gets:
\begin{align*}
\text{for a.e. $t \in \R \text{ such that } t\frac{x-P_1}{|x-P_1|} \in L_x \cap \omega$}, \quad 
\tilde{m}_{L_x}(t)= -n^{-1}_B \left( \frac{x-P_1}{|x-P_1|} \right), \\
\text{for a.e. $t \in \R \text{ such that } t\frac{z-P_1}{|z-P_1|} \in L_y \cap \omega$}, \quad 
\tilde{m}_{L_z}(t)= -n^{-1}_B \left( \frac{z-P_1}{|z-P_1|} \right),
\end{align*}
and using it to the lines $L_y$ and $L_z$, one gets:
\begin{align*}
\text{for a.e. $t \in \R \text{ such that } t\frac{y-P_2}{|y-P_2|} \in L_x \cap \omega$}, \quad 
\tilde{m}_{L_x}(t)= -n^{-1}_B \left( \frac{y-P_2}{|y-P_2|} \right), \\
\text{for a.e. $t \in \R \text{ such that } t\frac{z-P_2}{|z-P_2|} \in L_y \cap \omega$}, \quad 
\tilde{m}_{L_z}(t)= -n^{-1}_B \left( \frac{z-P_2}{|z-P_2|} \right).
\end{align*}
But $n^{-1}_B(-x)=-n^{-1}_B(x)$ and the two traces obtained for $\tilde{m}_{L_z}$ are opposed 
on the segment $[P_1,P_2]$, wich leads to a contradiction (see Figure \ref{figure:impossible}).

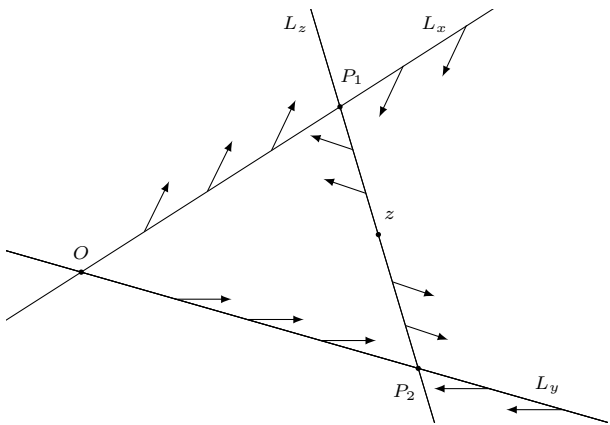
\begin{figure}[ht!]
\begin{tikzpicture}[line cap=round,line join=round,>=latex,x=1.0cm,y=1.0cm]
\clip(-1,-2) rectangle (7,3.5);
\draw [domain=-4.3:7.4] plot(\x,{(-0--2.2*\x)/3.44});
\draw [domain=-4.3:7.4] plot(\x,{(-0-1.28*\x)/4.48});
\draw [domain=-4.3:7.4] plot(\x,{(--14.26-3.48*\x)/1.04});
\draw [->] (0.84,0.54) -- (1.17,1.21);
\draw [->] (1.68,1.08) -- (2.01,1.75);
\draw [->] (2.53,1.62) -- (2.85,2.29);
\draw [domain=-4.3:7.4] plot(\x,{(-0-1.28*\x)/4.48});
\draw [domain=-4.3:7.4] plot(\x,{(-0-1.28*\x)/4.48});
\draw [->] (1.25,-0.36) -- (2,-0.36);
\draw [->] (2.21,-0.63) -- (2.96,-0.63);
\draw [->] (3.17,-0.91) -- (3.92,-0.91);
\draw [->] (4.28,2.74) -- (3.96,2.06);
\draw [->] (5.12,3.28) -- (4.8,2.6);
\draw [->] (5.44,-1.55) -- (4.69,-1.55);
\draw [->] (6.4,-1.83) -- (5.65,-1.83);
\draw [->] (4.31,-0.71) -- (4.88,-0.9);
\draw [domain=-4.3:7.4] plot(\x,{(--14.26-3.48*\x)/1.04});
\draw [->] (4.14,-0.13) -- (4.7,-0.32);
\draw [->] (3.61,1.63) -- (3.04,1.82);
\draw [->] (3.78,1.05) -- (3.22,1.24);
\begin{scriptsize}
\fill [color=black] (0,0) circle (1pt);
\draw[color=black] (0,0.26) node {$O$};
\fill [color=black] (3.44,2.2) circle (1pt);
\draw[color=black] (3.6,2.6) node {$P_1$};
\fill [color=black] (3.95,0.5) circle (1pt);
\draw[color=black] (4.1,0.76) node {$z$};
\fill [color=black] (4.48,-1.28) circle (1pt);
\draw[color=black] (4.3,-1.6) node {$P_2$};
\draw[color=black] (6.2,-1.5) node {$L_y$};
\draw[color=black] (4.7,3.3) node {$L_x$};
\draw[color=black] (2.85,3.3) node {$L_z$};
\end{scriptsize}
\end{tikzpicture}
\caption{An impossible configuration}
\label{figure:impossible}
\end{figure}

\textbf{Case 2}: $\text{dist}(O,\omega)>d$. Using the first part, we know that up to a change of sign:
$m(x)= n^{-1}_{B}\left( \frac{x}{|x|} \right)=V_B(x)$ and $m(y)= n^{-1}_{B}\left(\frac{y}{|y|}\right)=V_B(y)$, 
so that
\begin{align*}
|m(x)- m(y)|= \left|V_B(x)-V_B(y) \right|\leq \frac{K}{d^{\frac{1}{p-1}}} |x-y|^{\frac{1}{p-1}},
\end{align*}
by Proposition \ref{prop:convexity_regularity_body}.
This concludes our proof.
\end{proof}

\begin{rmk}
Notice that this result is sharp, in the sense that we cannot hope to get a better regularity for vector fields satisfying the kinetic equation and valued in the sphere of a power type $p$ norm. Indeed, as shown in Section \ref{sec:convexity} the vortex $V_B$ associated to the ball $B$ of this norm is always solution, and Proposition \ref{prop:convexity_regularity_body} shows that it is $\frac{1}{p-1}$-Hölder continuous far from $0$ if and only if the norm is of power type $p$.
\end{rmk}

\bigskip

\textbf{Acknowledgments}:
The first author would like to thank warmly Felix Otto for helpful discussions and comments
and the Max Planck Institute of Leipzig for the six months he spent there.
This work was supported by the LABEX MILYON (ANR-10-LABX-0070) of Université de Lyon, 
within the program "Investissements d'Avenir"
(ANR-11-IDEX-0007) operated by the French National Research Agency (ANR).

\bibliographystyle{plain}
\bibliography{Rot_nul_valeur_variete.bib}

\end{document}